\documentclass[12pt]{article}
\usepackage{setspace}
\usepackage[pdftex]{graphicx}
\usepackage{graphics, amsmath}
\usepackage{bm}
\usepackage[left=1.0in,top=1.0in,right=1.0in,bottom=1.0in]{geometry}

\usepackage{amssymb,amsmath,latexsym}

\def\independenT#1#2{\mathrel{\rlap{$#1#2$}\mkern2mu{#1#2}}}

\newcommand{\biblist}{\begin{list}{}
{\listparindent 0.0cm \leftmargin 0.50cm \itemindent -0.50 cm
\labelwidth 0 cm \labelsep 0.50 cm
\usecounter{list}}\clubpenalty4000\widowpenalty4000}
\newcommand{\ebiblist}{\end{list}}

\usepackage{latexsym}

\newcommand{\bx}{{\bf x}}
\newcommand{\bX}{{\bf X}}

\newcommand{\bB}{{\bf B}}
\newcommand{\bb}{{\bf b}}
\newcommand{\bD}{{\bf D}}
\newcommand{\bg}{{\bf g}}
\newcommand{\bG}{{\bf G}}

\newcommand{\bI}{{\bf I}}
\newcommand{\bW}{{\bf W}}
 
\newcommand{\btheta}{\mbox{\boldmath$\theta$}} 
\newcommand{\bbeta}{\mbox{\boldmath$\beta$}}

\newcommand{\bPhi}{\mbox{\boldmath$\Phi$}} 
\newcommand{\lk}{\left[} 
\newcommand{\rk}{\right]} 
\newcommand\independent{\protect\mathpalette{\protect\independenT}{\perp}}
\AtBeginDocument{
\addtolength{\abovedisplayskip}{-1ex}
\addtolength{\abovedisplayshortskip}{-1ex}
\addtolength{\belowdisplayskip}{-1ex}
\addtolength{\belowdisplayshortskip}{-1ex}
}
\newcommand{\be}{\begin{equation}} 
\newcommand{\ee}{\end{equation}} 

\begin{document}
\setcounter{page}{1}

\title{{\bf Generalized Method of Moments Estimator Based On Semiparametric Quantile Regression Imputation}}
\author{Senniang Chen\thanks{Department of Statistics, Iowa State University, Ames, IA 50011, USA. snchen@iastate.edu} and Cindy Yu\thanks{Department of Statistics, Iowa State University, Ames, IA 50011, USA. cindyyu@iastate.edu}}
\date{}
\maketitle

\begin{abstract}
In this article, we consider an imputation method to handle missing response values based on semiparametric quantile regression estimation. In the proposed method, the missing response values are generated using the estimated conditional quantile regression function at given values of covariates. We adopt the generalized method of moments for estimation of parameters defined through a general estimation equation. We demonstrate that the proposed estimator, which combines both semiparametric quantile regression imputation and generalized method of moments, has competitive edge against some of the most widely used parametric and non-parametric imputation estimators. The consistency and the asymptotic normality of our estimator are established and variance estimation is provided. Results from a limited simulation study and an empirical study are presented to show the adequacy of the proposed method.  
\end{abstract}

\noindent{\it Key Words}: generalized method of moments, imputation, semiparametric quantile regression.

\section{Introduction} Missing data is frequently encountered in many disciplines. Missing data analyses are important because inference based on ignoring missingness undermines efficiency and often leads to biases and misleading conclusions. The large literature handling missing data can basically be grouped into three categories: observed likelihood-based approaches, inverse probability weighting methods, and imputation methods. The main motivation for imputation is to provide a complete data set so that standard analytical techniques can be applied and the resulting point estimates are consistent among different users. Due to its intuitive simplicity, imputation becomes particularly popular among practitioners and is the focus of our paper. 

Many different imputation approaches have been developed in the literature and some prominent examples are included here. Rubin's pioneer work (1987) discussed multiple imputation (MI) based on Bayesian methods to generate pseudo values from the posterior predictive distribution and impute multiple data sets. Despite its simple form, however, the variance estimator of MI will have convergence problems if congeniality and self-sufficiency conditions are not met (Meng 1994). Fractional imputation was proposed to retain both estimation efficiency of multiple imputation and consistency of the Rao-Shao variance estimator (Rao and Shao 1992). In fractional imputation, multiple values are imputed for each missing cell with assigned weights. Kim (2011) proposed parametric fractional imputation (PFI) with inspirations from importance sampling and calibration weighting to reduce the computation burden. Noticeably, both PFI and MI assume a parametric regression model, and therefore may suffer from model misspecification. While MI and PFI resort to the creation of artificial responses, hot-deck imputation (HDI) replaces missing units with observed data through matching methods. By using covariate information, the matching method could be classifying donors and recipients into similar categorical classes (Brick and Kalton 1996; Kim and Fuller 2004), or creating metrics to match donors and recipients (Rubin 1986; Little 1988). More examples are documented in Andridge and Little (2010). In a recent work by Wang and Chen (2009), multiple imputed values are independently drawn from observed respondents with probabilities proportional to kernel distances between missing cells and donors. Both HDI and Wang and Chen (2009) are purely non-parametric, so the stability and accuracy of the estimators depend on the dimensionality and the sample size concerned. In fact, finite sample biases are observed in both of these non-parametric methods in our simulation study. It might be due to the fact that a donor with a higher probability of being present is more likely selected for imputing than a donor with a lower probability of being present, which possibly results in a distorted conditional density when the covariate is non-uniformly distributed. For more detailed discussions about this issue, see Section 3. 

To leverage the advantages of both parametric and non-parametric methods and avoid the limitation of a pure or exclusive approach, we propose an imputation method based on semiparametric quantile regression, which has the following set up. Define $f(y|\bf x)$ as the conditional density where $y$ is the response subject to missing and $\bf x$ is the covariate always observed, and  $q(\tau|\bf x)$ as the $\tau$-th conditional quantile function, which is the inverse conditional distribution function $F^{-1}(\tau|\bf x)$. Instead of estimating $f(y|\bf x)$ parametrically or non-parametrically, we estimate $q(\tau|\bf x)$ semiparametrically using observed data under the missing at random (MAR) assumption, in the sense intended by Rubin (1976). Then multiple imputed values $y^*_j (j=1,\cdots,J)$ are obtained via $y^*_j= \hat{q}(\tau_j|\bf x)$, where $\tau_j$ is independently drawn from Uniform$[0, 1]$. The semiparametric quantile regression imputation (hereafter called SQRI) is expected to have appealing features. Firstly, the entire conditional distribution function is used to draw imputed values, hence preserving the conditional density of the filled-in response values. Secondly, because different conditional quantiles instead of conditional means or actual observations are used in imputation, the method is  less sensitive to outliers, as quantiles are known to be less affected by extremes. Thirdly, it does not require strong model assumptions as in a fully parametric solution, and is therefore robust against model violations. Lastly, imputed values can be easily created through random numbers generated from Uniform$[0, 1]$ once $\hat{q}(\tau_j|\bf x)$ is estimated. 

In this paper, we are interested in estimating parameters defined through a general estimation equation. After imputation, the data set is regarded as complete and the generalized method of moments (GMM) is used for parameter estimation. However, combining GMM estimation with SQRI (hereafter called SQRI-GMM) has not been studied, to our best knowledge. So it is not clear, despite its aforementioned theoretical appeals, whether the proposed method can be advocated as an effective alternative in imputation. There are two main goals in this paper. The first goal is to rigorously establish a large sample theory of our GMM estimator based on SQRI, and the second goal is to evaluate its finite sample performance through numerical simulation. We examine the first goal in Section 2 and investigate the second goal in Section 3 through addressing the following three research questions: (1) Can our SQRI-GMM method significantly reduce biases caused by model misspecification, compared with MI and PFI? Our simulations are contrived to cover different kinds of misspecified mean structures, and performances of the estimators are compared; (2) Can our SQRI-GMM method have competitive finite sample performance, compared with some established non-parametric imputation methods? This question is interesting since both hot-deck imputation and Wang and Chen (2009) are also robust against model violations. (3) Can our SQRI-GMM method provide credible inference? The coverage probability of the confidence interval based on our SQRI-GMM estimator is studied in the simulation. Through the analyses of these three important questions, this paper demonstrates some numerical advantages of our SQRI-GMM estimator for imputation. 

We are not the first to use quantile regression for imputation. Papers pertaining to quantile regression imputation include Munoz and Rueda (2009), Wei et al (2012) and Yoon (2013). Our paper is distinctive from these papers in terms of objective, type of imputation and theory. (i) For objective, while Wei et al (2012) and Yoon (2013) limited their attention to the estimation of quantile regression coefficients, our method can be used for estimating parameters defined through any general estimation equation. Munoz and Rueda (2009) focused on the imputation strategy only, and parameter estimation was not an objective of the paper. It is worth noting that the setting in Wei et al (2012) is also different since they dealt with missing covariates, not missing responses. (ii) For type, Wei et al (2012) imputed multiple data sets,  while Munoz and Rueda (2009) proposed a single and deterministic imputation. However, our method utilizes fractional imputation. (iii) For theory, instead of assuming a linear quantile regression model, like Wei et al (2012) and Yoon (2013) did, we rely on a flexible semiparametric approach incorporating penalty for model complexity. And the key idea, which is used to arrive at the  consistency and normality in the proof, is substantially different from Wei et al (2012) and Yoon (2013). Because the primary interest of Munoz and Rueda (2009) was the computation strategy, no theory was offered in their paper. Our paper is unique in its contribution to theory building and its emphasis on application for a general framework with less restrictive assumptions.

The rest of paper is organized as follows. In Section 2, we introduce our imputation method through semiparametric quantile regression with penalty and present large sample theories of our SQRI-GMM estimator. Section 3 compares our method with some competing methods through simulation studies and reports the statistical inference results of our SQRI-GMM estimator.  Section 4 analyzes an income data set from Canadian Census Public Use Tape. The Appendix outlines proofs of the theorems appearing in the main text. Details of the proofs are in the supplemental file Chen and Yu (2014).

\section{Proposed GMM Estimator Based On Semiparametric Quantile Regression Imputation (SQRI)} In this section, we introduce our GMM estimator based on SQRI. Section 2.1 builds the framework and discusses SQRI using penalized B-splines. Section 2.2 establishes the asymptotic consistency and normality of the unweighted SQRI-GMM estimator. Section 2.3 extends the large sample theory to the weighted SQRI-GMM estimator.

\subsection{SQRI using penalized B-splines} We consider $(\bx_i, y_i)^T, i=1,\cdots,n$, to be a set of i.i.d. observations of random variables $(\bX, Y)$, where $\bf{X}$ is a $d_x$-dimension variable always observed and $Y$ is the response variable subject to missing. Let $\delta_i=1$ if $y_i$ is observed and $\delta_i=0$ if $y_i$ is missing. We assume that $\delta$ and $Y$ are conditionally independent given $\bf{X}$, i.e. 
$$
P(\delta=1|Y=y, \bX=\bx)=P(\delta|\bX=\bx) := p(\bx),
$$
a condition called ``missing at random'' by Rubin (1976). The primary interest of this article is to estimate a $d_\theta$-dimensional parameter $\btheta_0$, which is the unique solution to $E\{\bg(Y,\bX; \btheta)\}=0$, and make inference on $\btheta_0$. Here $\bg(Y,\bX; \btheta)=(g_1(Y,\bX; \btheta),\cdots,g_r(Y,\bX; \btheta))^T$ is a vector of $r$ estimating functions for $r\geq d_{\theta}$. Let $q_{\tau}(\bx)$ be the unknown conditional $100\tau$\% quantile of response $Y$ given $\bX=\bx$. It satisfies
$
P(Y<q_{\tau}(x)|\bX=\bx)=\tau
$
for a given $\tau \in (0,1)$. When $\tau=0.5$, $q_{\tau}(\bx)$ is the conditional median of $Y$. It is easy to show that $q_{\tau}(\bx)$ satisfies
$$
q_{\tau}(\bx) =\arg\min_{h(\bx)} E\{\rho_{\tau}(Y-h(\bx))|\bX=\bx \},
$$
where $\rho_{\tau}(u)=u(\tau-I(u<0))$, the check function proposed in Koenker and Bassett (1978). Many have studied the estimation of $q_{\tau}(\bx)$ based on parametric methods, and a summary of relevant literature can be found in Koenker (2005). Parametric model assumptions may not hold sometimes, giving rise to nonparametric methods. Nonparametric quantile regression, including the kernel quantile regression in Yu and Jones (1994) and the smoothing spline method in Koenker et al (1994), has also been intensively studied. Among many findings is the well-known trade-off between computational cost and smoothness. In other words, spline smoothing methods demand massive computation, and the unpenalized spline tends to give wiggly curves despite its cheap computational cost.  In this paper, we employ a semiparametric quantile regression method based on penalized B-splines, as suggested in Yoshida (2013), that features a relatively smoothed quantile function at reduced computational burden. 

To simplify notations, we assume $X$ is an univariate variable with a distribution function $F_x(x)$ on $[0, 1]$. We discuss how to deal with multivariate $\bX$ in Section 3 and how to rescale $\bX$ on any compact set into $[0, 1]$ in Section 4. Let $K_n -1$ be the number of knots within the range $(0, 1)$, and $p$ be the degree of B-splines. In order to construct the $p$-th degree B-spline basis, we define equidistantly located knots as $\kappa_k= K_n^{-1} k, (k=-p+1,\cdots,K_n+p)$. Note there are $K_n-1$ knots located in $(0,1)$. The $p$-th B-spline basis is 
$$
\bB(x)=(B^{[p]}_{-p+1}(x), B^{[p]}_{-p}(x),\cdots,B^{[p]}_{K_n}(x))^T,
$$  
where $B^{[p]}_{k}(x) (k=-p+1,\cdots,K_n)$ are defined recursively as 
\begin{itemize}
\item For $s=0$:
$$
B^{[s]}_{k}(x)=B^{[0]}_{k}(x)= \left \{
\begin{array}{ll}
1, \kappa_{k-1}<x \leq \kappa_k,\\
0, \mbox{ otherwise},\\
\end{array}
\right. 
\mbox{ where } k=-p+1,\cdots,K_n+p.
$$
\item For $s=1,2,\cdots,p$:
$$B_{k}^{[s]}(x)= \frac{x-\kappa_{k-1}}{\kappa_{k+s-1}-\kappa_{k-1}}B_{k}^{[s-1]}(x)+ \frac{\kappa_{k+s}-x}{\kappa_{k+s}-\kappa_{k}}B_{k+1}^{[s-1]}(x),$$ 
\hspace{0.4in} where $ k=-p+1,...,K_n+p-s.$
\end{itemize}

Readers can refer to de Boor (2001) for more details and properties of the B-spline functions. The estimated conditional quantile regression function is $\hat{q}_{\tau}(x)=\bB^T(x) \hat{\bb}(\tau)$, where $\hat{\bb}(\tau)$ is a $(K_n+p) \times 1$ vector obtained by 
\be
\label{b-hat}
\hat{\bb}(\tau)= \arg\min_{\bb(\tau)} \sum_{i=1}^n \delta_i \rho_\tau[y_i-\bB^T(x_i) \bb(\tau)]+\frac{\lambda_n}{2}\bb^T(\tau) \bD_m^T \bD_m \bb(\tau). 
\ee
Here $\lambda_n (>0)$ is the smoothing parameter, and $\bD_m$ is the $m$-th difference matrix and is $(K_n+p-m) \times (K_n+p)$ dimensional with its element defined as 
$$d_{ij}=\left\lbrace\begin{array}{ll}(-1)^{|i-j|}\left(\begin{array}{c}m\\|i-j|\end{array} \right) & 0\leq j-i\leq m\\
 0& o.w. \end{array}\right. ,
$$
where the notation $\left(\begin{array}{c} m\\k\end{array} \right)$ is the choose function given by $(k! (m-k)!)^{-1} m!$ and $m$ is the order of penalty. As discussed in Yoshida (2013), the difference penalty $\bb^T(\tau) \bD_m^T \bD_m \bb(\tau)$ is used to remove computational difficulty occurring when the penalty term is defined through an integral, and it controls the smoothness of the estimated quantile regression function. Section 3 discusses how we choose the numbers $(\lambda_{n},m,K_n,p)$ in practice. 

To control the variability of the estimating functions with imputed values, we generate $J$ independent imputed values $\{ y_{ij}^*\}_{j=1}^J$ when $y_i$ is missing by the following procedure.
\begin{enumerate}
\item Simulate $\tau_j \sim$ Uniform(0,1) independently for $j=1, 2,\cdots, J$;
\item For each $j=1, 2,..., J$, $\hat{\bb}(\tau_j)$ is calculated as 
$$\hat \bb(\tau_j)= \arg\min_{\bb(\tau)} \sum_{i=1}^n \delta_i \rho_{\tau_j}[y_i- \bB^T(x_i) \bb(\tau)]+\frac{\lambda_n}{2}\bb^T(\tau) \bD_m^T \bD_m \bb(\tau);
$$
\item For the missing unit $i$, $J$ independent values are generated as
$$y^*_{ij}| x_i=\hat{q}_{\tau_j}(x_i)=\bB^T(x_i)\hat{\bb}(\tau_j),j=1, 2,\cdots, J.$$
\end{enumerate}
Repeat step 3 for every missing unit in the data set. Then we use 
$
\delta_i \bg(x_i,y_i; \btheta)+(1-\delta_i) J^{-1} \sum_{j=1}^J \bg(x_i,y^*_{ij}; \btheta)
$
as the estimating function for the $i$-th observation. 

Sometimes the conditional mean of $Y$ given $X=x$ is used for imputation, such as in Cheng (1994) and Wang and Rao (2002), but it does not work for general parameter estimation. For some  parametric imputation methods, imputation and estimation steps are entwined, in that updating parameters and re-imputing based on most recently updated parameters are iteratively done. This might require heavy computing time. In the SQRI described above, imputation and estimation steps are totally separate, making general purpose parameter estimation possible. Also in SQRI, standard analytical tools can be directly applied to imputed data without re-imputation. The PFI by Kim (2011) avoids re-imputation by adjusting weights of imputed values based on iteratively updated parameters. However, any parametric imputation method, including PFI and MI, might suffer from model misspecification. Non-parametric imputation, such as HDI or the method proposed in Wang and Chen (2009) using kernel distance, assumes no parametric model, but the stability and accuracy of non-parametric estimators depend on sample size and dimensionality of the problem. The SQRI provides a useful compromise between a fully parametric approach and a purely non-parametric approach. 

Assume the number of knots $K_n-1$ and the smoothing parameter $\lambda_{n}$ depend on $n$. By Barrow and Smith (1978), there exists $\bb^*(\tau)$ that satisfies 
\be \label{barrow-smith}
\sup_{x \in (0,1)} | q_\tau(x)+b_{\tau}^a(x)-\bB^T(x) \bb^*(\tau) |=o(K_n^{-(p+1)}),
\ee
where $b_\tau^a(x)= \frac{q_\tau^{(p+1)}(x)}{(p+1)!K_n^{p+1}}Br_p(\frac{x-\kappa_{k-1}}{K_n^{-1}})$ if $\kappa_{k-1}\leq x < \kappa_k$, and $q_\tau^{(p+1)}(x)$ is the $(p+1)$-th derivative of $q_\tau(x)$ with respect to $x$. Here $Br_p(\cdot)$ is the $p$-th Bernoulli polynomial inductively defined as
$Br_0(x)=1, \text{ and }  Br_p(x)=\int_0^x p B_{p-1}(z)dz+b_p,$
where $b_p=-p\int_0^1\int_0^xBr_{p-1}(z)dzdx$ is the $p$-th Bernoulli number (Barrow and Smith (1978) and Yoshida (2013)). The following Lemma gives the asymptotic property of $\hat{q}_{\tau}(x)=\bB^T(x) \hat{\bb}(\tau)$ where $\hat{\bb}(\tau)$ is defined in (\ref{b-hat}). 

\vspace{0.2in}

\noindent \textbf{Lemma 1: } Under condition 1 given in the Appendix, and assuming $q_{\tau}(x) \in C^{p+1}$,  $K_n=O(n^{\frac{1}{2p+3}})$, and $\lambda_n =O(n^v)$ for $v \leq (2p+3)^{-1}(p+m+1)$, we have  
\\
(i) 
\be \label{Lemma1-i}
\sqrt{\frac{n}{K_n}}[\hat q_\tau(x)-\bB^T(x)\bb^*(\tau)+b^\lambda_\tau(x)]\to_d N(0,V_\tau), 
\ee
(ii)
\be \label{Lemma1-ii}
\sqrt{\frac{n}{K_n}}[\hat{q}_\tau(x)-q_\tau(x)+b^a_\tau(x)+b^\lambda_\tau(x)]\to_d N(0,V_\tau), 
\ee
for a given $x \in (0,1)$ and $\tau \in (0,1)$, where  
\be
\begin{array}{lll}
b^\lambda_\tau(x)&=& \frac{\lambda_n}{n}\bB^T(x)\left(\bPhi(\tau)+\frac{\lambda_n}{n}\bD^T_m \bD_m\right)^{-1}\bD^T_m \bD_m \bb^*(\tau),\\
 V_\tau(x) &=& \lim_{n \to \infty}\frac{\tau(1-\tau)}{K_n}\bB^T(x)\left(\bPhi(\tau)+\frac{\lambda_n}{n}\bD^T_m \bD_m\right)^{-1}\bPhi\\
 &&\times \left(\bPhi(\tau)+\frac{\lambda_n}{n}\bD^T_m \bD_m\right)^{-1}\bB(x),\\
\bPhi&=&\int_0^1 p(x) \bB(x) \bB^T(x) dF_x(x),\\
\bPhi(\tau)&=&\int_0^1 p(x) f_{y|x}(q_{\tau}(x)) \bB(x) \bB^T(x) dF_x(x).\\
\end{array}
\ee
Here $f_{y|x}(\cdot)$ is the conditional density of $Y$ given $X=x$. There exist two sources of asymptotic biases in $\hat{q}_\tau(x)$. One is $b^a_\tau(x)$ which is the model bias between the true function $q_{\tau}(x)$ and the spline model used, see equation (\ref{barrow-smith}). Another source of bias $b^\lambda_\tau(x)$ is introduced by adding the penalty term into the quantile regression. When there is no penalty term ($\lambda_n=0$), this bias vanishes. Both of these two bias terms have an order $O(n^{-\frac{p+1}{2p+3}})$. The proof of this lemma draws from Theorem 1 of Yoshida (2013), which deals with complete data. The detailed proof of this order and Lemma 1 can be found in the supplemental file Chen and Yu (2014).  

We define $\bG(\btheta)=E\{\bg(Y,X; \btheta)\}$ and our estimating function as
\be
\bG_n(\btheta)=\frac{1}{n} \sum_{i=1}^n \{\delta_i \bg(y_i,x_i;\btheta)+(1-\delta_i) \frac{1}{J} \sum_{j=1}^J \bg(y_{ij}^*, x_i; \btheta) \}.
\ee
We consider the generalized method of moments (GMM), a usual estimation equation approach, to make inference on $\btheta$. Our proposal of combining SQRI with GMM is attractive, thanks to the applicability to general parameter estimation of GMM and the aforementioned appeals of SQRI. 

\subsection{Unweighted GMM estimator based on SQRI} The unweighted GMM-SQRI estimator is obtained as
\be \label{unweighted}
\widehat{\btheta}_n=\arg\min_{\btheta \in \Theta } \bG_n^T(\btheta) \bG_n(\btheta).
\ee
We first present Lemma 2, which regards the asymptotic normality of $\bG_n(\btheta_0)$.

\vspace{0.2in}

\noindent \textbf{Lemma 2: } Under conditions 1 and 2 (a) $\sim$ (b) given in the Appendix, and assuming $q_{\tau}(x) \in C^{p+1}$,  $K_n=O(n^{\frac{1}{2p+3}})$, and $\lambda_\tau =O(n^v)$ for $v \leq (2p+3)^{-1}(p+m+1)$, as $n\rightarrow \infty$ and $J\rightarrow \infty$ we have  
\be
\sqrt{n}\bG_n(\btheta_0) \to_d N(0,V_G(\btheta_0)), 
\ee
where 
\be \label{V-G}
V_G(\btheta)=Var(\xi_i(\btheta)),
\ee
\be \label{xi}
\xi_i(\btheta)= \bg(y_i,x_i;\btheta)+(1-\delta_i)\lk \mu_{g|x}(x_i; \btheta)-\bg(y_i,x_i;\btheta) \rk+\delta_iC_p h_n(y_i,x_i; \btheta) \bB(x_i),
\ee
\be \label{hn}
h_n(y_i,x_i;\btheta)=\int_{-\infty}^{+\infty} \int_{0}^1  \dot \bg_y(q_{\tau}(x),x; \btheta) \bB^T(x) H_n^{-1}(\tau) \psi_{\tau}(e_i(\tau)) d\tau dF_X(x),
\ee
$$
\mu_{g|x}(x; \btheta)=E\{\bg(y,x; \btheta)|X=x \}, \mbox { } H_n(\tau)=\Phi(\tau)+\frac{\lambda_n}{n}\bD_m^T\bD_m, 
$$
$$
\dot\bg_y(y,x; \btheta)=\frac{\partial \bg(y,x; \btheta)}{\partial y}, \mbox{   } \psi_{\tau}(u)=\tau-1_{u<0},
$$
$$
e_i(\tau)=y_i-\bB^T(x_i) \bb^*(\tau),  \mbox{ and }  C_p=E\{ 1-p(x)\}.
$$
Justification of Lemma 2 is crucial to show consistency and asymptotic normality of our SQRI-GMM estimator (Pakes and Pollard 1989). We decompose $\sqrt{n}\bG_n(\btheta_0)$ into three terms
\be \label{Bterms}
\begin{array} {lll}
\sqrt{n}\bG_n(\btheta_0)& = &\underbrace{\frac{1 }{\sqrt n}\sum_{i=1}^n \bg(y_i,x_i;\btheta_0)}_{:=B_1}\\
&&+ \underbrace{ \frac{1 }{\sqrt n} \sum_{i=1}^n[(1-\delta_i) ( \mu _{g|x}(x_i;\btheta_0)-\bg(y_i,x_i;\btheta_0)]}_{:=B_2}\\
                        &   & +\underbrace{\frac{1 }{\sqrt n} \sum_{i=1}^n[(1-\delta_i)(\hat \mu_{g|x}(x_i;\btheta_0)-\mu_{g|x}(x_i;\btheta_0)]}_{:=B_3},\\
\end{array}
\ee
where $  \hat \mu_{g|x}(x_i;\btheta) = J^{-1}\sum_{j=1}^J \bg(y_{ij}^*,x_i, \btheta)$ and $y_{ij}^*=\bB^T(x_i) \hat{\bb}(\tau_j)$. The terms $B_1$ and $B_2$ are simple since they are sums of i.i.d. random variables. The term $B_3$ is much more complicated because it involves additional randomness from the uniformly distributed random variable $\tau_j$, and it also depends on the estimated coefficients $\hat{\bb}(\tau_j)$ calculated using all respondents. Therefore the summands in $B_3$ are not independent. The key idea in the proof is to replace $B_3$ by $\tilde{B}_3=E(B_3|A_R)$ where $A_R=\{\delta_i, (y_i,x_i)|\delta_i=1; i=1,\cdots,n  \}$, and to show the following two results:
 (1)
  $
\tilde{B}_3=n^{-1/2} \sum_{i=1}^n \delta_i C_p h_n(y_i,x_i; \btheta_0) \bB(x_i)+o_p(1)
$, and (2) $\tilde{B}_3-B_3=o_p(1)$. Combining these two results with equation (\ref{Bterms}) gives the asymptotic normality in Lemma 2. 

\vspace{0.1in}
 
\noindent \textit{Remark 1:} When there is no missing, $\xi_i(\btheta_0)$ in equation (\ref{xi}) coincides with $\bg(y_i,x_i;\btheta_0)$. 

\vspace{0.2in}

\noindent \textbf{Theorem 1: } Under conditions 1 and 2 given in the Appendix, and assuming $q_{\tau}(x) \in C^{p+1}$,  $K_n=O(n^{\frac{1}{2p+3}})$, and $\lambda_\tau =O(n^v)$ for $v \leq (2p+3)^{-1}(p+m+1)$, as $n\rightarrow \infty$ and $J\rightarrow \infty$ we have  

(i) $$ \widehat{\btheta}_n \to_p \btheta_0; $$

(ii) $$\sqrt{n} \Sigma^{-1/2}(\btheta_0) (\widehat{\btheta}_n-\btheta_0) \to_d N(\bf{0}, \bI_{d_{\theta} \times d_{\theta} }), $$

where 
\be \label{Sigma}
\Sigma(\btheta)= \left\{\Gamma^T(\btheta)\Gamma(\btheta)\right\}^{-1}\Gamma^T(\btheta) V_G(\btheta) \Gamma(\btheta)\left\{\Gamma^T(\btheta)\Gamma(\btheta)\right\}^{-1},\ee
\be \mbox{ and } \Gamma(\btheta)=E\{\frac{\partial \bg(Y,X; \btheta)}{\partial \btheta}\}.\ee

Theorem 1 shows that $\widehat{\btheta}_n$ is consistent and asymptotically normal. With Lemma 2 and the fulfillment of the following 2 conditions: (1) $\sup_{\btheta} (1+|\bG(\btheta)|+|\bG_n(\btheta)|)^{-1} |\bG_n(\btheta)-\bG(\btheta)|=o_p(1)$ and (2) $\sup_{|\btheta-\btheta_0|<\zeta_n} (n^{-1/2} + |\bG(\btheta)|+|\bG_n(\btheta)|)^{-1}|\bG_n(\btheta)-\bG(\btheta)-\bG_n(\btheta_0)|=o_p(1)$ for any positive sequence $\zeta_n$ converging to zero, Theorem 1 can be proved following Corollary 3.2 and Theorem 3.3 of Pakes and Pollard (1989). Here the notation of $|\cdot|$ represents the norm of a matrix, defined as $|A|=\sqrt{trace(A^{\prime} A)}$. 

To consider variance estimation for $\widehat{\btheta}_n$, let an estimator of $\xi_i(\btheta)$ be
\be
\hat{\xi_i}(\btheta)= \bg(y_i,x_i;\btheta)+(1-\delta_i)\left\{\hat{\mu}_{g|x}(x_i; \btheta)-\bg(y_i,x_i;\btheta)\right\}+\delta_i \hat C_p \hat{h}_n(y_i,x_i; \btheta) \bB(x_i),
\ee
where   
$$ \hat{h}_n(y_i,x_i;\btheta) = \frac{1}{n}\frac{1}{J}\sum_{k=1}^n \sum_{j=1}^J \dot \bg_y(\hat q_{\tau_j}(x_k),x_k; \btheta) \bB^T(x_k) \hat H_n^{-1}(\tau_j)\psi_{\tau_j}(\hat e_i(\tau_j)),$$
 $$ \hat e_i(\tau_j)  = y_i-\bB^T(x_i) \hat \bb(\tau_j) \mbox{, }
 \hat H_n(\tau_j) =  \hat \Phi(\tau_j)+\frac{\lambda_n}{n}\bD_m^T\bD_m, $$
 $$\hat \Phi(\tau_j) = \frac{1}{n}\sum_{i=1}^n \delta_i \hat f_{Y|X}(x_i,\hat q_{\tau_j}(x_i)) \bB(x_i)\bB^T(x_i) \mbox{ with } \hat q_{\tau_j}(x_i)=\bB^T(x_i)\hat{\bb}(\tau_j),$$
 $$\hat f_{Y|X}(x,y) = \frac{\frac{1}{nab}\sum_{i=1}^n \delta_i K(\frac{y-y_i}{a})K(\frac{x-x_i}{b})}{\frac{1}{na}\sum_{i=1}^n \delta_i K(\frac{x-x_i}{a})}, \mbox{ and }\hat C_p = n^{-1}\sum_{i=1}^n(1-\delta_i).$$
Here the estimation of $\hat f_{Y|X}(x,y)$ uses a Normal kernel $K(\cdot)$ and bandwidth $a$ or $b$ for $x$ (or $y$). The estimator of $\Gamma(\btheta_0)$ is obtained by
$$\hat{\Gamma}(\hat \btheta_n)=\frac{1}{n}\sum_{i=1}^n \left\{\delta_i \frac{\partial \bg(y_i,x_i;\hat \btheta_n)}{\partial \btheta} +(1-\delta_i) \frac{1}{J}\sum_{j=1}^J  \frac{\partial \bg(y^*_{ij},x_i;\hat \btheta_n)}{\partial \btheta}\right\}.$$
Then, the variance estimator of $\widehat{\btheta}_n$ is $\hat{V}(\widehat{\btheta}_n )=n^{-1} \widehat{\Sigma}(\widehat{\btheta}_n)$, where $\widehat{\Sigma}(\widehat{\btheta}_n)$ is calculated as
\be
\widehat{\Sigma}(\widehat{\btheta}_n)=\left\{\hat\Gamma^T(\widehat{\btheta}_n) \hat \Gamma(\widehat{\btheta}_n)\right\}^{-1}\hat \Gamma^T(\widehat{\btheta}_n) \hat{V}_G(\widehat{\btheta}_n) \hat \Gamma(\widehat{\btheta}_n)\left\{\hat \Gamma^T(\widehat{\btheta}_n) \hat \Gamma(\widehat{\btheta}_n)\right\}^{-1}, \mbox{ and }
\ee
\be \label{Vhat-G}
 \hat{V}_G(\btheta) = \frac{1}{n-1}\sum_{i=1}^n \left\{\hat \xi_i(\btheta)-\frac{1}{n}\sum_{i=1}^n\hat \xi_i(\btheta)\right\}\left\{\hat \xi_i(\btheta)-\frac{1}{n}\sum_{i=1}^n\hat \xi_i(\btheta)\right\}^T.
\ee

\vspace{0.2in}

\noindent \textbf{Corollary 1: } Under conditions 1 $\sim$ 3 given in the Appendix, and assuming $q_{\tau}(x) \in C^{p+1}$,  $K_n=O(n^{\frac{1}{2p+3}})$, and $\lambda_\tau =O(n^v)$ for $v \leq (2p+3)^{-1}(p+m+1)$, as $n\rightarrow \infty$ and $J\rightarrow \infty$ we have  
$$\sqrt{n} \widehat{\Sigma}^{-1/2}(\widehat{\btheta}_n) (\widehat{\btheta}_n-\btheta_0) \to_d N(\bf{0}, \bI_{d_{\theta} \times d_{\theta} }). $$
Corollary 1 allows us to construct confidence intervals based on the asymptotic normality and the variance estimator.

\subsection{Weighted GMM estimator based on SQRI} A weighted GMM estimator is calculated by minimizing $\bG_n^T(\btheta) \bW \bG_n(\btheta)$ for a $r \times r$ positive definite weight matrix $\bW$. It can be shown that taking $\bW \propto V^{-1}_G(\btheta_0)$ will result in the most efficient estimator in the class of all asymptotically normal estimators using arbitrary weight matrices. In practice, $\bW$ is replaced by the inverse of the random matrix $ \hat{V}_G(\btheta)$ defined in (\ref{Vhat-G}) and the weighted GMM estimator is obtained as 
\be \label{weighted}
\widehat{\btheta}^w_n=\arg\min_{\btheta \in \Theta } \bG_n^T(\btheta) \hat{V}^{-1}_G(\btheta) \bG_n(\btheta).
\ee
The following Lemma proves that $\hat{V}^{-1}_G(\btheta) $ is close to the fixed non-singular matrix $V^{-1}_G(\btheta_0)$ uniformly over a sequence of shrinking neighborhoods, an important condition for $\widehat{\btheta}^w_n$ to be consistent and asymptotically normal. 

\vspace{0.2in}

\noindent \textbf{Lemma 3: } Under conditions 1 $\sim$ 3  given in the Appendix, and assuming $q_{\tau}(x) \in C^{p+1}$,  $K_n=O(n^{\frac{1}{2p+3}})$, and $\lambda_\tau =O(n^v)$ for $v \leq (2p+3)^{-1}(p+m+1)$, as $n\rightarrow \infty$ and $J\rightarrow \infty$ we have  
$$\sup_{|\btheta-\btheta_0|<\zeta_n} 
| \hat{V}^{-1}_G(\btheta)- V^{-1}_G(\btheta_0)|=o_p(1),$$
for a sequence of positive numbers $\zeta_n$ that converges to zero. 

The following theorem presents the large sample properties of the weighted GMM estimator $\widehat{\btheta}^w_n$. 

\vspace{0.2in}

\noindent \textbf{Theorem 2: } Under conditions 1 $\sim$ 3 given in the Appendix, and assuming $q_{\tau}(x) \in C^{p+1}$,  $K_n=O(n^{\frac{1}{2p+3}})$, and $\lambda_\tau =O(n^v)$ for $v \leq (2p+3)^{-1}(p+m+1)$, as $n\rightarrow \infty$ and $J\rightarrow \infty$ we have  

(i) $$ \widehat{\btheta}^w_n \to_p \btheta_0; $$

(ii) $$\sqrt{n} \Sigma_w^{-1/2}(\btheta_0) (\widehat{\btheta}^w_n-\btheta_0) \to_d N(\bf{0}, \bI_{d_{\theta} \times d_{\theta} }), $$

where $ \label{Sigma}\Sigma_w(\btheta)= \left\{\Gamma^T(\btheta) V_G^{-1}(\btheta) \Gamma(\btheta)\right\}^{-1}.$

When Lemma 3 holds, the results in Theorem 2 follow immediately from Lemmas 3.4 and 3.5 of Pakes and Pollard (1989).
 
\vspace{0.1in}
 
\noindent \textit{Remark 2:} The asymptotic variance of the most efficient GMM estimator based on the complete data is $n^{-1}[\Gamma^T(\btheta_0) Var^{-1}\{\bg(y_i,x_i;\btheta_0)\} \Gamma(\btheta_0)]^{-1} $. It can be shown that $V_G(\btheta_0)$ in equation (\ref{V-G}) can also be expressed as
\be \begin{array} {lll}
V_G(\btheta_0)&=& Var\{\bg(y,x;\btheta_0)\}-E\{(1-p(x))\sigma^2_{g|x}(x; \btheta_0)\}\\
&& +C_p^2E\left\{\delta_i h_n(y_i,x_i;\btheta_0)\bB(x_i) \bB^T(x_i)h_n^T(y_i,x_i;\btheta_0) \right\}\\
&&+2C_pE\left\{\delta_i h_n(y_i,x_i; \btheta_0)\bB(x_i) \bg^T(y_i,x_i;\btheta_0) \right\},
 \end{array}
\ee 
where $\sigma^2_{g|x}(x; \btheta_0)=Var\{\bg(y,x;\btheta_0)|X=x\}$. So when missing is low, i.e. $p(x)$ is large and $C_p$ is close to zero, the efficiency of $\widehat{\btheta}^w_n$ is close to the asymptotic efficiency of the best GMM estimator under no missing. 

\vspace{0.1in}

\noindent \textit{Remark 3:} When $r=d_{\theta}$, the semiparametric efficiency bound defined in Chen, Hong and Tarozzi (2008) is $$\Sigma_{speb}(\btheta_0)=[\Gamma^T(\btheta_0) E^{-1}\{\sigma^2_{g|x}(x;\btheta_0)/p(x)+\mu_{g|x}(x;\btheta_0) \mu^T_{g|x}(x;\btheta_0) \}\Gamma(\btheta_0)]^{-1}.$$ Rewrite 
\begin{eqnarray*}
V_{G}(\btheta_0)&=&E\{p(x) \sigma^2_{g|x}(x;\btheta_0)\}+V\{\mu_{g|x}(x,\btheta_0)\}\\
&&+C_p^2E\{\delta_i h_n(y_i,x_i;\btheta_0)\bB(x_i) \bB^T(x_i)h_n^T(y_i,x_i;\btheta_0)\}\\
&&+2C_pE\{\delta_i h_n(y_i, x_i; \btheta_0)\bB(x_i) \bg^T(y_i,x_i;\btheta_0)\}.
\end{eqnarray*} Our estimator will achieve the semiparametric efficiency bound if $V_G(\btheta_0) \leq E\{\sigma^2_{g|x}(x;\btheta_0)/p(x)+\mu_{g|x}(x;\btheta_0) \mu^T_{g|x}(x;\btheta_0) \}$, i.e. 
\be \label{criteria}
\begin{array}{lll}
E\{(\frac{1}{p(x)}-p(x))\sigma^2_{g|x}(x;\btheta_0)\}&\geq& C_p^2E\{\delta_i h_n(y_i,x_i;\btheta_0)\bB(x_i) \bB^T(x_i)h_n^T(y_i, x_i;\btheta_0)\}\\
&& +2C_pE\{\delta_i h_n(y_i, x_i;\btheta_0)\bB(x_i) \bg^T(y_i,x_i;\btheta_0)\}.
\end{array}
\ee
It can be shown that under the conditions $K_n=O(n^{\frac{1}{2p+3}})$ and $\lambda_\tau =O(n^v)$ for $v \leq (2p+3)^{-1}(p+m+1)$, the right hand side of equation (\ref{criteria}) has order $O(K_n^{-1})$ (see derivation in Chen and Yu (2014)). However, the left side is $O(1)$. So when $K_n \rightarrow \infty$, inequality (\ref{criteria}) will likely happen. This might explain why in our simulation studies our estimator has slightly smaller Monte Carlo variances than the non-parametric imputation estimator of Wang and Chen (2009), which is claimed to have the semiparametric efficiency bound when $r=d_{\theta}$.

The variance estimator for $\widehat{\btheta}^w_n$ can be simply computed as $\hat{V}(\widehat{\btheta}^w_n )=n^{-1}\widehat{\Sigma}_w(\widehat{\btheta}^w_n)$, where $\widehat{\Sigma}_w(\widehat{\btheta}^w_n)= \{\hat{\Gamma}(\widehat{\btheta}^w_n)^T \hat{V}^{-1}_G(\widehat{\btheta}^w_n) \hat{\Gamma}(\widehat{\btheta}^w_n)\}^{-1}.$ The following Corollary shows that the central limit theory still holds after replacing $\Sigma_w(\btheta)$ by its estimator, thus  inference can be legitimately made based on the weighted SQRI-GMM estimator and its variance estimator. 

\vspace{0.2in}

\noindent \textbf{Corollary 2: } Under conditions 1 $\sim$ 3  given in the Appendix, and assuming $q_{\tau}(x) \in C^{p+1}$,  $K_n=O(n^{\frac{1}{2p+3}})$, and $\lambda_\tau =O(n^v)$ for $v \leq (2p+3)^{-1}(p+m+1)$, as $n\rightarrow \infty$ and $J\rightarrow \infty$ we have  
$$\sqrt{n} \widehat{\Sigma}_w^{-1/2}(\widehat{\btheta}^w_n) (\widehat{\btheta}^w_n-\btheta_0) \to_d N(\bf{0}, \bI_{d_{\theta} \times d_{\theta} }). $$

\section{Simulation Studies} The second goal of our article is to evaluate the finite sample performances of our proposed estimator through simulation studies. For this purpose, we investigate the following three research questions: (i) Can our proposed method significantly reduce  biases caused by model misspecification, compared with parametric imputation methods such as MI and PFI? (ii) How does our proposed method perform, compared with non-parametric imputation methods such as hot-deck imputation and the method proposed in Wang and Chen (2009)? (iii) Can credible inference be made based on our proposed method? 

We specify the simulation set-up as follows. The response $y_i$ is generated from a model $y_i=m(\bx_i)+\epsilon_i$, where $m(\bx_i)$ is the mean function and $\epsilon_i$ are $iid$ $N(0,0.1^2)$. We consider the following four different mean functions drawing from the design of simulation studies in Breidt et al (2005) to cover a range of correct and incorrect model specification. 
$$
\begin{array} {rl}
\mbox{\textit{linear}: } & m(x_i)=1+2(x_i-0.5),\\
\mbox{\textit{bump}: } & m(x_i)=1+2(x_i-0.5)+\exp\{-30(x_i-0.5)^2 \},\\
\mbox{\textit{cycle}: } & m(x_i)=0.5+2x_i+sin(3\pi x_i),\\
\mbox{\textit{bivariate}: } & m(x_{1i}, x_{2i})=1+2(x_{1i}-0.5)+2\exp\{-10(x_{2i}-0.4)^2 \}.\\
\end{array}
$$
The covariate $x_i$ for the first three univariate models (or $x_{1i}$ and $x_{2i}$ for the last \textit{bivariate} model) are all independently and identically simulated from a truncated normal distribution $N(0.5,0.3^2)$ on interval $[0,1]$. The missing mechanism considered follows a logistic regression model
$$
\begin{array} {lrlll}
&p(x_i)&=&\frac{\exp(1+0.5x_i)}{1+\exp(1+0.5x_i)} & \mbox{ for the \textit{linear, bump, cycle}  models},\\
\mbox{ or } &p(x_{1i},x_{2i})&=& \frac{\exp(0.2+x_1+0.5x_2)}{1+\exp(0.2+x_1+0.5x_2 )} & \mbox{ for the  \textit{bivariate} model}.\\
\end{array}
$$
The missing rates in all situations are about 20\%. We are interested in estimating three parameters, the marginal mean of the response variable $\mu_y=E(Y)$, the marginal standard deviation of the response variable $\sigma_y=\sqrt{Var(Y)}$ and the correlation between the response and covariate variables $\rho=corr(X,Y)$. So $\btheta=(\mu_y, \sigma_y, \rho)$ and the corresponding estimating function is defined as
\begin{equation} \label{define g}
\bg(x_i,y_i,\mu_x,\mu_y,\sigma_x,\sigma_y,\rho)=\left( \begin{array}{c}  x_i-\mu_x \\
                                                     y_i-\mu_y\\
                                                     (x_i-\mu_x)^2-\sigma_x^2\\
                                                     (y_i-\mu_y)^2-\sigma_y^2\\
                                                     (x_i-\mu_x)(y_i-\mu_y)-\rho \sigma_x\sigma_y\\  \end{array} \right).
\end{equation}
For  \textit{bivariate} model, $\btheta=(\mu_y, \sigma_y, \rho_1, \rho_2)$, where $\rho_1=corr(X_1,Y)$ and $\rho_2=corr(X_2,Y)$ and the estimating function is defined in an analogous way. Note that $\mu_x$ and $\sigma^2_x$ are the mean and variance of the covariate and are treated as nuisance parameters. If there is no missing, the parameter vector $\btheta$ can be estimated as
\begin{equation} \label{no-missing}
\begin{array} {lll} 
\hat{\mu}_y=\frac{1}{n} \sum_{i=1}^n y_i, & \hat{\sigma}^2_y=\frac{1}{n-1} \sum_{i=1}^n (y_i-\hat{\mu}_y)^2, &\\
\hat{\mu}_x=\frac{1}{n} \sum_{i=1}^n x_i, & \hat{\sigma}^2_x=\frac{1}{n-1} \sum_{i=1}^n (x_i-\hat{\mu}_y)^2,& \\
\hat{\rho}=\frac{n^{-1} \sum_{i=1}^n (x_i-\hat{\mu}_x)(y_i-\hat{\mu}_y)}{\hat{\sigma}_x \hat{\sigma}_y}. && \\
\end{array}
\end{equation}
For each model, 1000 replicate samples of size $n=200$ are created and the following seven estimators are calculated to compare our semiparametric imputation method with some established parametric and non-parametric imputation methods. 
\begin{itemize}
\item \textbf{Full}: An estimator based on the full observations. $\hat{\btheta}$ is calculated using equation (\ref{no-missing}).
\item \textbf{Resp}: A naive estimator based on respondents only (where ``Resp'' comes from the word ``respondents''). $\hat{\btheta}$ is calculated using equation (\ref{no-missing}) after ignoring missing. 
\item \textbf{SQRI-GMM}: Our proposed estimator defined in (\ref{weighted}), which combines the semiparametric quantile regression imputation and weighted GMM estimation. 
\item \textbf{MI}: The multiple imputation estimator proposed in Rubin (1987). The $R$ package `mi' by Gelman et al (2013) is employed to obtain $J$ multiple imputed data sets. Estimators in (\ref{no-missing}) are calculated for each imputed data set, and the MI estimators are averaged across multiple imputed data sets. 
\item \textbf{PFI}: The parametric fractional imputation estimator proposed in Kim (2011). Under PFI, multiple imputed values $y_{ij}^* (j=1,\cdots,J)$ are generated from a proposed conditional density $\tilde{f}(y|x)$ and their associated fractional weights $w_{ij}^*$ are computed using $\tilde{f}(y|x)$ and the assumed conditional density $f(y|x; \hat{\eta}^0)$, where $\hat{\eta}^0$ is the given initial value for $\hat{\eta}$ in the conditional density formula. By maximizing the score function of the density $f(y_i|x_i; \eta)$ using the imputed values and their weights, $\hat{\eta}$ is updated, and the fractional weights $w_{ij}^*$ are re-calculated iteratively until $\hat{\eta}$ converges. The PFI estimators are calculated using equation (\ref{no-missing}), with the missing $y_i$ values replaced by $\sum_{j=1}^J w_{ij}^* y_{ij}^*$.
\item \textbf{NPI-EL}: The non-parametric imputation estimator based on the empirical likelihood method in Wang and Chen (2009). In NPI-EL, multiple imputed values $y_{ij}^* (j=1,\cdots,J)$ are independently drawn from the respondent group $(\delta_i=1)$ with the probability of selecting $y_s$ with $\delta_s=1$ being
$$
P(y_{ij}^*=y_s)=\frac{K\{(x_s-x_i)/h\}}{\sum_{m=1}^n \delta_m K\{(x_m-x_i)/h\}},
$$
where $K(\cdot)$ is a $d_x$-dimensional kernel function and $h$ is a smoothing bandwidth. In our simulations, the Gaussian kernel is used with $h$ prescribed by a cross-validation method. The NPI-EL estimator is obtained using the empirical likelihood method for a general estimation problem with the estimating function for a missing unit $i$ replaced by $J^{-1} \sum_{j=1}^J \bg(y_{ij}^*, x_i;\btheta)$.
\item \textbf{HDFI}: A hot-deck fractional imputation estimator. Under HDFI, multiple imputed values $y_{ij}^* (j=1,\cdots,J)$ are independently drawn from a donor pool consisting of 20 nearest neighbors identified through the Euclidean distance in the current study. The HDFI estimators are calculated using (\ref{no-missing}) with the missing $y_i$ replaced by $J^{-1} \sum_{j=1}^J y_{ij}^*$. 
\end{itemize}

The Full and the Resp estimators are included in order to help us gauge how far away our proposed estimator is from the ideal case and from the case of simply ignoring missing. Estimators NPI-EL and HDFI are based on non-parametric imputation methods, while estimators MI and PFI are based on parametric imputation methods, where $y_i$ is assumed to satisfy $Y|X=x \sim N(\bbeta^T \bx, \sigma^2)$ for some $\sigma>0$. Our SQRI-GMM is semiparametric as we use penalized B-spline to estimate conditional quantile regression. For penalized B-spline quantile estimators, typically the degree of B-spline $p$ and the degree of the difference matrix $m$ are fixed at low values, for example $p \leq 3$ and $m\leq 2$. We set $p=3$ and $m=2$, a popular choice in practice as suggested in Yoshida (2013). For a given $K_n$ (where $K_n=\mbox{ number of knots}+1$), the smoothing parameter $\lambda_{n}$ is prescribed via the generalized approximation cross-validation (GACV) method discussed by Yuan (2006). We obtain results for a variety of choices of $K_n$ and conclude $K_n=5$ suffices in our examples. In the \textit{bivariate }model, the same specifications are used to obtain bases $\bB(x_1)$ and $\bB(x_2)$ on $x_1$ and $x_2$ separately, then $\bB(\bx)$ is their augmentation, $\bB(\bx)=(\bB^T(x_1), \bB^T(x_2))^T.$ For all the five imputation methods described above, we use both $J=10$ and $J=100$. The simulation studies show that $J=10$ is sufficient for our proposed estimator to accurately estimate parameters.  We summarize the numerical findings for $J=10$ below. Conclusions are the same for $J=100$.

Table 1-2 present the Monte Carlo relative biases and variances of the seven estimators for the four models. To compare bias, we compute the ratios of relative biases for other estimators and the relative bias for the proposed SQRI-GMM estimator, and take the absolute values. If the absolute ratio is bigger than 1, the proposed estimator has smaller relative biases. Figure 1 is the visualization of the bias comparison. The relative biases of the proposed estimator are less than 1\% in all cases and are closest to those of the Full estimator in nearly all cases (Table 1 and 2). In particular, the proposed estimator has smaller biases and variances as compared with the Resp estimator because the former incorporates additional covariate information of the missing units while the latter totally ignores missing units. 

The following findings are summarized to answer research Question (i), which addresses the performance of the proposed estimator with respect to some parametric imputation estimators. When the \textit{linear} model is correctly specified, the SQRI-GMM estimator has relative biases of a magnitude comparable to the two parametric estimators MI and PFI. When the model is misspecified (\textit{bump, cycle, bivariate}), our estimator has significantly smaller biases than the MI and PFI estimators with exceptions arising in $\mu_y$ and $\sigma_y$ of the \textit{bivariate} model, where  the three estimators all have small relative biases less than 1\%. Correspondingly, Figure 1 (b) $\sim$ (d) show the proposed estimator to be advantageous with relative ratios mostly over the threshold of 1 and reaching as high as 50, compared with MI (curve with squares) and PFI (curve with triangles). In terms of variances (Table 1 and 2), the proposed semi-parametric estimator has slightly bigger variances than the two parametric ones under the correct \textit{linear} model as expected, but has slightly smaller or close efficiency under the incorrect models.  

The following findings are summarized to answer research Question (ii), which addresses the performance of the proposed estimator with respect to some non-parametric imputation estimators. Compared with the two non-parametric estimators NPI-EL and HDFI, our estimator has considerably smaller biases with only one exception when estimating $\rho_2$ in the \textit{bivariate} model, where our relative bias is -0.0070 and the relative bias of the NPI-EL is 0.0056. This superior performance can be seen in Figure 1(a) $\sim$ (d) where the curves with circles (NPI-EL) and the curves with stars (HDFI) are well above the horizontal line of 1 reaching as high as 80. The variances of the proposed estimator are generally in line with that of the HDFI estimator (Table 1 and 2). When compared with the NPI-EL estimator, the SQRI-GMM estimator has slightly smaller varainces. This corresponds to Remark 3 in Section 2.3 where the condition is $K_n \rightarrow \infty$. However, this superiority in efficiency is small because we only use $K_n=5$ in practice. 

The biases observed in the two non-parametric methods can be possibly explained by the fact that respondents with a higher probability of being present are more likely selected for imputing than respondents with a lower probability of being present when $x$ is non-uniformly distributed. An artificial example is plotted in Figure 2 to help with illustration. This example mimics the \textit{linear} model used in the simulation where the covariate $x$ follows a truncated normal distribution centered at $x=0.50$, and the units with higher $x$ values have higher probabilities of being present. Suppose we want to impute the $y$ value at $x=0.25$ using HDFI and assume there is no observation between $x\in (0.12,0.38)$, an illustrative situation to facilitate explanation. The donor group consists of 10 nearest neighbors (highlighted bigger dots) that are at the same distance away from $x=0.25$. There are 9 respondents around $x=0.40$ and only 1 respondent at $x=0.10$ due to the  non-uniform distribution of $x$. The location of $J^{-1} \sum_{j=1}^J y_{ij}^*$ at $x=0.25$ calculated from the 10 donors is marked by the symbol $*$ in Figure 2. These imputed values will pull the true conditional mean up, resulting in overestimation of $\mu_y$. It is consistent with the findings in Table 1 and 2 that both NPI-EL and HDFI overestimate the marginal mean $\mu_y$ in all cases. Similar overestimating effect will occur if there are observations between $x\in (0.12,0.38)$ because there will be more donors on the right side of $x=0.25$ than on the left side of $x=0.25$ for the same reason. This argument can also explain the biases associated with NPI-EL. Under the NPI-EL method, the 10 highlighted dots have the same chance of being drawn as imputed value because they have the same kernel distances away from $x=0.25$. Therefore, more imputed values will be from those 9 respondents at $x=0.40$, resulting in a bigger $J^{-1} \sum_{j=1}^J y_{ij}^*$ value. In fact, Table 1 and 2 show that NPI-EL and HDFI have the same directions of over or under estimation across models and parameters. Another possible reason is that both NPI-EL and HDFI are arguably local methods which might occasionally suffer from unstable estimates in regions with high missing rates. However, our estimator is based on global quantile regression, and thus is less sensitive to the presence of such regions relative to purely non-parametric methods. 

The following findings are summarized to answer Question (iii), which is about the inference validity of the proposed estimator. Table 3 contains the coverage probabilities of our SQRI-GMM estimator based on asymptotic normality (Corollary 2) and a bootstrapping method for both $J=10$ and $J=100$. For $J=10$, the coverage probabilities based on normality are close to the nominal level of 0.95 except  for $\rho$ under the \textit{linear} and \textit{cycle} models. This is common for confidence intervals constructed based on normal approximation of a GMM estimator. After increasing from $J=10$ to $J=100$, all coverage probabilities based on normality improve in general, though the coverages for $\rho$ in the \textit{linear} and \textit{cycle} models are still low (about 86\% and 92\%). A bootstrapping method then is conducted as a remedy to obatain the confidence intervals. The bootstrapping algorithm is described as follows. 

\begin{enumerate}
  \item Draw a simple random sample $\chi^*_n$  with replacement from the original sample $\chi_n= {(X_i,Y_i,\delta_i)_{i=1}^n}$;
 \item Implement semiparametric quantile regression to impute values for the missing cells in $\chi^*_n$;
 \item Estimate $\hat\btheta$ using the SQRI-GMM estimator.
 \item Repeat step 1 $\sim $ 3 for B times, then we have $\hat\btheta^1, \hat\btheta^2, \cdots, \hat\btheta^B$.
 \end{enumerate} 
The $2.5$-th and $97.5$-th percentiles of $\{\hat\btheta^b\}_{b=1}^B$ give the lower and upper bounds of the $95\%$ confidence interval. We use $B=400$ in our simulation. In general, the bootstrapping method has a slightly better performance over normal approximation method, offering satisfactory coverage probabilities close to 0.95 even when $J$ is small. 

In summary, our simulation studies confirm the validity of our proposed estimator in  finite sample estimation.  

\section{Empirical Study} In this section, our proposed SQRI-GMM estimator is applied to a real data set consisting of $n=205$ Canadian workers all educated to grade 13. A description of the data set can be found in Ruppert et al (2003) and Ullah(1985), by whom the source was identified as a 1971 Canadian Census Public Use Tape. A copy of the data can be obtained from the $R$ package `SemiPar' by Wand (2013). The study variable $y$ is the natural logarithm of \textit{annual income} and the covariate $x$ is \textit{age} rescaled into $[0,1]$ by the formula $x=(age-min(age))/(max(age)-min(age))$. The sample estimates of $(\mu_y,\sigma_y,\rho)$  are $(13.49, 0.636, 0.231)$ when there is no missing.  Missingness is created artificially by deliberately deleting some of the $y$ values according to the missing mechanism $p(x)=exp(1-0.5x)/\{1+ exp(1-0.5x)\}$, which results in a 30\% missing rate. All the five imputation estimators described in the simulation are computed using the real data with artificial missing.

The variance estimator for MI is a function of the point estimators and the variance estimators based on all imputed data sets. We use GMM to obtain both point and variance estimators for each imputed data set. The variance estimators for PFI and HDFI are computed using a bootstrapping method similar to what was described in Section 3 except that different imputation methods are employed in Step 2. The confidence interval for NPI-EL is obtained via the bootstrapping method introduced in Wang and Chen (2009). Table \ref{table4} reports the relative biases (relative to the sample estimates of $(\mu_y,\sigma_y,\rho)$ based on full observations) and 95\% confidence interval widths for five estimators. Figure 4 is the scatterplot of $income$ on a log scale versus $age$. When estimating $\mu_y$, all estimators have relative biases less than 1\%. However, when estimating $\sigma_y$ and $\rho$, there exists telling differences: the relative biases of our estimator are smaller than those of other estimators. This might be due to some features of the data. For example, there is no obvious mean structure (pattern) after age 22, which might explain why all estimators can estimate the overall mean well; also there is noticeable heteroscedasticity in the data, which might cause other estimators to fall short.  In general, our estimator has slightly narrower confidence intervals except when estimating $\rho$ (it is inferior to the MI estimator). Overall, this case study demonstrates the empirical effectiveness of the SQRI-GMM estimator.

\section*{Acknowledgements}
The authors thank Cooperative Agreement No. 68-3A75-4-122 between the USDA Natural Resources Conservation Service and the Center for Survey Statistics and Methodology at Iowa State University.

\section*{Appendix} The notation of $|\cdot|$ represents the norm of a matrix, defined as $|A|=\sqrt{trace(A^{\prime} A)}$ and the notation of $\left\|\cdot \right\|$ denotes the sup-norm in all arguments for functions. We first discuss some technical assumptions. 

1. Assumptions for penalized semiparametric quantile regression:
(a) There exists $\gamma>0$ such that $E[|g(y,x;\btheta)|^{2+\gamma}] < \infty$ .
(b) The explanatory variable $X$ has distribution function $F_x(x)$ on a compact set $[0,1]$.
(c) The knots for the B-spline basis are equidistantly located as $\kappa_k=k/K_n$ for $k=-p+1, \cdots, K_n+p$.
(d) The order of the difference matrix  is $m < p$.
(e) $\lim_{n \to \infty}n^{-1}\sum_{i=1}^np(x_i)B(x_i)B^T(x_i)$ exists and converges to $\Phi$ where $\Phi$ is defined as
$\Phi=\int_0^1\bB(x) \bB^T(x)p(x)dF_x(x).$
(f) $\lim_{n \to \infty} n^{-1}\sum_{i=1}^n p(x_i) f_{Y|X}(q_\tau(x))$ $\bB(x_i)\bB^T(x_i) $ exists and converges to $\Phi(\tau)$, where $\Phi(\tau)=\int_0^1\bB(x) \bB^T(x)p(x)$ $f_{Y|X}(q_\tau(x))dF_x(x).$
(g) The smoothing parameters $\lambda_n$ is a positive sequence of real numbers such that $\lambda_n^{-1}$ is larger than the maximum eigenvalue of $\Phi(\tau)^{-1/2}D_m^TD_m\Phi(\tau)^{-1/2}$. 

2. Assumptions for the GMM method:
(a) $\btheta_0$ is the unique solution to the general estimating equation $E[\bg(x,y,\btheta)]=0$ and $\btheta_0$ is an interior point of $\Theta$.
(b) $\bg(y,x,\btheta)$ is differentiable with respect to $\btheta$ and twice differentiable with respect to $y$. $\dot \bg_\theta(y,x,\btheta)=\frac {\partial \bg(y,x,\btheta) }{\partial \btheta}$, $\ddot \bg_{\theta,y}(y,x,\btheta)=\frac {\partial^2 \bg(y,x,\btheta) }{\partial \btheta\partial y}$ and $\ddot \bg_y(y,x,\theta)=\frac {\partial^2 \bg(y,x,\btheta_0) }{\partial y^2}$ are bounded for all $\btheta \in \Theta$, $x$, and $y$. 
(c) $\bG(\btheta)=E\lk \bg(y,x;\btheta)\rk$ is differentiable at $\btheta_0$ with a derivative matrix $\Gamma$ of full rank.
(d) $||\bG_n(\hat\btheta_n)|| \leq o_p(n^{-1/2})+\inf_{\btheta}|\bG_n(\btheta)|$.
(e) $ E\lk |\dot \bg_\theta(y,x;\btheta)||\dot \bg_\theta(y,x;\btheta)|^T\rk$ is bounded.

3. General assumptions for variance estimators:
(a) The bandwidths of the kernel density estimator, $a$ and $b$, satisfy $a \to 0$, $b \to 0$, $na \to \infty$ and $nb \to \infty$.
(b) $f_{Y|X}(x,y)$ is differentiable with respect to y.
(c) $\| \Gamma(\btheta)\|$ and $\| V_G(\btheta)\|$ are bounded away from 0, where $V_G(\btheta)=Var\lk \xi_i(\btheta)\rk$ and $\xi_i(\btheta)= \bg(y_i,x_i;\btheta)+(1-\delta_i)\lk \mu_{g|x}(x_i,\btheta)-\bg(y_i,x_i;\btheta)\rk+\delta_i C_p h_n(x_i,y_i,\btheta)\bB(x_i)$.

The proof for Lemma 1 and Corollary 2 are skipped here because the proof for Lemma 1 is very similar to Theorem 1 of Yoshida (2013) except that we are dealing with missing data and the proof of Corollary 2 is straightforward and similar to the proof for Corollary 1. We give outlines of the proofs for the rest of the theories stated in the text in this Appendix. More detailed proofs (including the skipped ones), the facts referred hereafter and their justifications can be found in the supplemental file Chen and Yu (2014). 
\subsection*{A: Proof of Lemma 2}
\renewcommand{\theequation}{A.\arabic{equation}}
\renewcommand{\thesection}{A}
\setcounter{equation}{0} 
We can decompose $\bG_n(\btheta_0)$ as in (\ref{Bterms}). The key idea in our proof is to replace $B_3$ by $\tilde{B}_3=E(B_3|A_R)$ where $A_R=\{\delta_i, (y_i,x_i)|\delta_i=1;i=1,\cdots,n  \}$, and to show the following two results: (1)
$
\tilde{B}_3=n^{-1/2} \sum_{i=1}^n\delta_i C_p h_n(x_i,y_i,\btheta_0)\bB(x_i)+o_p(1),$
and (2)
$\tilde{B}_3-B_3=o_p(1)$. 

\noindent \textbf{(1)} To show $\tilde B_3=n^{-1/2} \sum_{i=1}^n\delta_i C_p h_n(x_i,y_i,\btheta_0)\bB(x_i)+o_p(1)$: We  further decompose $\tilde B_3$ into two terms,
\be
\begin{array}{lll}
\tilde B_3&=&\underbrace{\frac{1 }{\sqrt n}\sum_{i=1}^n E\left\lbrace (1-\delta_i)\frac{1}{J}\sum_{j=1}^J[\bg(\hat q_{\tau_j}(x_i),x_i;\btheta_0)-\bg(q_{\tau_j}(x_i),x_i;\btheta_0)]|A_R  \right\rbrace }_{\tilde B_{31}}\\
&& +\underbrace{\frac{1 }{\sqrt n}\sum_{i=1}^nE\left\lbrace (1-\delta_i)\frac{1}{J}\sum_{j=1}^J[\bg(q_{\tau_j}(x_i),x_i;\btheta_0)-\mu_{g|x}(x_i;\btheta_0)]|A_R  \right\rbrace }_{\tilde B_{32}}.
\end{array}
\ee
It is obvious that $\tilde B_{32}=0$ because $E_{\tau|x}\lk \bg(q_{\tau}(x),x,\btheta_0)\rk=E_{y|x}\lk \bg(y,x,\btheta_0) \rk =\mu_{g|x}(x;\btheta_0)$ for any $x$.
For $\tilde B_{31}$, assuming that $\bg(x_i,y_i,\btheta)$ is twice differentiable with respect to $\btheta$, then we have 
\be\label{gtalyor}\begin{array}{lll}
&&\bg(\hat q_{\tau}(x_i),x_i;\btheta_0)-\bg(q_{\tau}(x_i),x_i;\btheta_0)\\&= &\dot \bg_y(q_{\tau}(x_i),x_i;\btheta_0)[\hat q_{\tau}(x_i)- q_{\tau}(x_i)] + \ddot \bg_y(\tilde q_{\tau}(x_i),x_i;\btheta_0)[\hat q_{\tau}(x_i)- q_{\tau}(x_i)]^2, 
\end{array}
\ee
 for $\tilde q_{\tau_j}(x_i)$ lying between $q_{\tau_j}(x_i)$ and $\hat q_{\tau_j}(x_i)$. By equation (\ref{gtalyor}), we have
\begin{eqnarray*}
\tilde B_{31}&=&\frac{n_m}{\sqrt n} E\left\lbrace \frac{1}{J}\sum_{j=1}^J\dot \bg_y(q_{\tau_j}(x),x,\btheta_0)[\hat q_{\tau_j}(x)-q_{\tau_j}(x)]|A_R  \right\rbrace \\
&& +  \frac{n_m}{\sqrt n} E\left\lbrace \frac{1}{J}\sum_{j=1}^J\ddot \bg_y(\tilde q_{\tau_j}(x),x;\btheta_0)[\hat q_{\tau_j}(x)-q_{\tau_j}(x)]^2|A_R \right\rbrace,
\end{eqnarray*}
where $n_m=n-\sum_{i=1}^n\delta_i$ and $x \independent A_R$. By Fact 3 in Chen and Yu (2014), we have
\be
\label{ddiff g}
E\left\lbrace \frac{1}{J}\sum_{j=1}^J\ddot \bg_y(\tilde q_{\tau_j}(x),x,\btheta_0)[\hat q_{\tau_j}(x)-q_{\tau_j}(x)]^2|A_R \right\rbrace =O(\frac{K_n}{n}).
\ee
By Lemma 1, we have
\be 
\label{qhat-q}\begin{array}{lll}
\sqrt n\left( \hat q_\tau(x)- q_\tau(x) \right)&=&\frac{1}{\sqrt n}\bB^T(x)H_n^{-1}(\tau)\sum_{i=1}^n\delta_i\bB(x_i)\psi_\tau(e_i(\tau))\\&&-\frac{\lambda_n}{\sqrt n}\bB^(x)H_n^{-1}(\tau)C_n(\tau)-\sqrt n b^a_\tau(x)+o_p(1),\end{array}
\ee
where $C_n(\tau)=\bD_m^T\bD_m\bb^*(\tau)$. Then $\tilde B_{31}$ can be written as
\be
\tilde B_{31}= \frac{C_p}{\sqrt n}\sum_{i=1}^n\delta_i h_n(x_i,y_i,\btheta_0)\bB(x_i)-\sqrt n C_pC_{1n}- \sqrt n C_pC_{2n}+o_p(1),
\ee
where $h_n(x_i,y_i,\btheta_0)$, $C_{1n}$, $C_{2n}$ and $C_p$ are defined in Lemma 2. The asymptotic order of $C_{1n}$ and $C_{2n}$ are 
$ C_{1n}=\frac{\lambda_n}{n} E_{x,\tau}\left\lbrace \dot \bg_y(q_{\tau}(x),x,\btheta)\bB^T(x)H_n^{-1}(\tau)C_n(\tau) \right\rbrace = O(K_n^{-(p+2)})$ by Fact 2, and $C_{2n}= E[\dot \bg_y(q_{\tau}(x),x,\btheta_0)\bb_\tau^a(x)]=O(K_n^{-(p+2)}).$ Thus we have
$
\tilde B_{31}=\frac{1}{\sqrt n}\sum_{i=1}^n\delta_i C_p h_n(x_i,y_i,\btheta_0)\bB(x_i)+o_p(1).
$

\noindent \textbf{(2)} To show $\tilde B_3-B_3=o_p(1)$: By Chebychev's inequality, we only need to show that
$E[\tilde B_3-B_3]^{\otimes2} \leq  E\{[\hat \mu_{g|x}(x_i;\btheta_0)-\mu_{g|x}(x_i;\btheta_0)]^{\otimes2} \}=o(1).
$
First of all, we can decompose $\hat \mu_{g|x}(x_i;\btheta_0)-\mu_{g|x}(x_i;\btheta_0)$ into two terms,
\begin{eqnarray*}
&&\hat \mu_{g|x}(x_i;\btheta_0)-\mu_{g|x}(x_i;\btheta_0)\\
&=&\underbrace{ \frac{1}{J}\sum_{j=1}^J[\bg(\hat q_{\tau_j}(x),x;\btheta_0)- \bg(q_{\tau_j}(x),x;\btheta_0)]}_{S_n}+\underbrace{ \frac{1}{J}\sum_{j=1}^J[\bg(q_{\tau_j}(x),x;\btheta_0)- \mu_{g|x}(x_i;\btheta_0)]}_{Q_n}.\end{eqnarray*}
It is equivalent to show that $E[Q_n^{\otimes2}]=o(1)$, $E[Q_nS^T_n]=o(1)$ and $E[S_n^{\otimes2}]=o(1)$. Details to show these orders can be found in Chen and Yu (2014). Combining step (1) \& 2, together with equation (\ref{Bterms}), we can write $\sqrt n\bG_n(\btheta_0)=\frac{1}{\sqrt n}\sum_i^n \xi_i(\btheta) + o_p(\frac{1}{\sqrt n}) $ where $\xi_i(\btheta)$ is defined in (\ref{xi}). Then by the central limit theorem, we have $V^{-1/2}(\xi_i(\btheta_0))n^{-1/2} \sum_{i=1}^n \xi_i(\btheta_0) \sim_d N(\bf{0}, \bI_{r \times r }),$
where $ V(\xi_i(\btheta_0))$ $\doteq \sigma_g^2(\btheta_0)-E[(1-p(x))\sigma_{g|x}^2(x;\btheta_0)]+C_p^2E\{p(x_i) h_n(x_i,y_i,\btheta_0)\bB(x_i) \bB^T(x_i)$ $h_n^T(x_i,y_i,\btheta_0) \} + 2C_pE\{\delta_i h_n(x_i,y_i)\bB(x_i) \bg^T(y_i,x_i;\btheta_0) \}$.

\subsection*{B: Proof of Theorem 1}
\renewcommand{\theequation}{B.\arabic{equation}}
\renewcommand{\thesection}{B}
\setcounter{equation}{0} 
 We verify the two conditions stated after Theorem 1. From equation (\ref{Bterms}) and $\tilde B_3-B_3=o_p(1)$ in Lemma 2, we have
$$\sqrt n\left(\bG_n(\btheta)-\bG(\btheta)\right)=B_1(\btheta)+ B_2(\btheta)+\tilde B_3(\btheta)+o_p(1),$$
where $B_2(\btheta)=n^{-1/2}\sum_{i=1}^n \lk (1-\delta_i)(\mu_{g|x}(x_i,\btheta)-g(y_i,x_i;\btheta)\rk$, $B_3(\btheta) = n^{-1/2}$  $\sum_{i=1}^n (1-\delta_i)(\hat\mu_{g|x}(x_i,\btheta)-\mu_{g|x}(x_i,\btheta)$ and $\tilde B_3(\btheta)=n^{-1/2} \sum_{i=1}^n\delta_iC_ph_n(x_i,y_i,\btheta)\bB (x_i)$. By the law of large numbers, we have $\frac{1}{n}\sum_i^n \{ \bg(y_i,x_i;\btheta)-E[(\bg(y_i,x_i;\btheta)]\}=o_p(1)$,  $\frac{1}{n}\sum_i^n \{ (1-\delta_i) ( \mu _{g|x}(x_i;\theta)-\bg(y_i,x_i;\btheta)\}=o_p(1)$ and $\frac{1}{n}\sum_i^n  \delta_i C_p  h_n(x_i,y_i;\btheta)$ $\bB(x_i) = o_p(1)$. Thus we have
$\| \bG_n(\btheta)-\bG(\btheta)\|=o_p(1)$ and $\sup_{\btheta}\frac{|\bG_n(\btheta)-\bG(\btheta) |}{1+|\bG_n(\btheta)|+|\bG(\btheta)|} $ $\leq \|\bG_n(\btheta)-\bG(\btheta) \| =o_p(1).$
So condition (1) holds. 

To prove condition (2), it is sufficient to show that for every sequence $\{\zeta_n\}$ of positive numbers converging to zero, $\bG_n(\btheta)-\bG(\btheta)-\bG_n(\btheta_0)=o_p(n^{-1/2}) \text{ for } \|\btheta-\btheta_0\|<\zeta_n.$ Since $\bG(\btheta_0)=0$, then
\be\label{GMMCond2}
       \begin{array}{lll}
      && \bG_n(\btheta)-\bG(\btheta)-\bG_n(\btheta_0)\\&=& \frac{1}{n} \sum_i^n [\bg(y_i,x_i;\btheta)-E(\bg(y,x;\btheta)]+B_2(\btheta)+ B_3(\btheta)\\
       &-&\left\{ \frac{1}{n}\sum_i^n [\bg(y_i,x_i;\btheta_0)-E(\bg(y,x;\btheta_0)]+B_2(\btheta_0)+ B_3(\btheta_0)\right\}+o_p(\frac{1}{\sqrt n}).
       \end{array}
 \ee
Because $|\btheta-\btheta_0|<\zeta_n$, we can show (details in Chen and Yu (2014)) that
$$\label{GMMCond21}  
   \frac{1}{n} \sum_i^n [\bg(y_i,x_i;\btheta)-E(\bg(y,x;\btheta)]-\frac{1}{n}\sum_i^n [\bg(y_i,x_i;\btheta_0)-E(\bg(y,x;\btheta_0)]= o_p(\frac{1}{\sqrt n}),$$
$      
       \frac{1}{\sqrt n} B_2(\btheta)-\frac{1}{\sqrt n}B_2(\btheta_0)= o_p(\frac{1}{\sqrt n}),
$
$
   \mbox{ and } \frac{1}{\sqrt n}\tilde B_3(\btheta)-\frac{1}{\sqrt n}\tilde B_3(\btheta_0)=o_p(\frac{1}{\sqrt n}).
$
Thus, we have $|| \bG_n(\btheta)-\bG(\btheta)-\bG_n(\btheta_0)||=o_p(\frac{1}{\sqrt n}),$ for $\{\zeta_n\}\rightarrow 0.$ 
                
\subsection*{C: Proof of Corollary 1}
\renewcommand{\theequation}{C.\arabic{equation}}
\renewcommand{\thesection}{C}
\setcounter{equation}{0} 
It is sufficient to show
$\hat V_G(\hat\btheta_n)=\hat Var(\hat\xi_i(\hat\btheta_n))$ converges to $V_G(\btheta_0)= Var(\xi_i(\btheta_0))$, $\hat \Gamma(\hat\btheta_n) \to^p\Gamma(\btheta_0) \text{ and }\left[\hat\Gamma^T(\hat \btheta_n)\hat\Gamma(\hat \btheta_n)\right]^{-1} \to_p \left[\Gamma^T( \btheta_0)\Gamma( \btheta_0)\right]^{-1}.
$
Here $\hat Var\left(\hat\xi(\hat\btheta_n)\right)=\frac{1}{n-1}\sum_{i=1}^n\left(\hat \xi_i(\hat\btheta_n)-\frac{1}{n}\sum_{i=1}^n\hat \xi_i(\hat\btheta_n)\right)^{\otimes2}$, which is a consistent estimator of $Var\left(\xi_i(\btheta_0)\right)$, if $\hat\xi_i(\hat\btheta_n) \to^p \xi_i(\btheta_0)$. The rest of the proof is to show that $\hat\xi_i(\hat\btheta_n) \to^p \xi_i(\btheta_0)$ by showing that $\bg(y_i,x_i,\hat\btheta_n)\to^p \bg(y_i,x_i,\btheta_0)$, $\hat\mu_{g|x}(x_i,\hat\btheta_n)\to^p \mu_{g|x}(x_i,\btheta_0) $ and $\hat h_n(x_i,y_i;\hat\btheta_n)\to^p   h_n(x_i,y_i;\btheta_0)$. Similarly, because of $\hat q_\tau(x) \to^p q_\tau(x)$ and $\hat\btheta_n \to^p \btheta_0$, as $n\to \infty$, and $J\to \infty$, we have
$\hat\Gamma(\hat \btheta_n) \to^p \Gamma(\btheta_0),$ and  when $\Gamma(\btheta_0)$ is bounded away from $\boldmath 0_{d_\theta\times d_\theta}$ for all $\btheta \in \Theta$, we have
$\left[\hat\Gamma^T(\hat \btheta_n)\hat\Gamma(\hat \btheta_n)\right]^{-1} \to^p \left[\Gamma^T( \btheta_0)\Gamma( \btheta_0)\right]^{-1}. $ Details can be found in Chen and Yu (2014). 

\subsection*{D: Proof of Lemma 3}  
\renewcommand{\theequation}{D.\arabic{equation}}
\renewcommand{\thesection}{D}
\setcounter{equation}{0} 
Since
$\sup_{| \btheta-\btheta_0| \leq \zeta_n}  | \hat V^{-1}_G(\btheta) - V^{-1}_G(\btheta_0)|  \leq   \sup_{| \btheta-\btheta_0| \leq \zeta_n}$ $ | \hat V^{-1}_G(\btheta)-V^{-1}_G(\btheta)| + \sup_{| \btheta-\btheta_0| \leq \zeta_n} | V^{-1}_G(\btheta)-V^{-1}_G(\btheta_0)|,
$ 
we only need to show that both $ \sup_{| \btheta-\btheta_0| \leq \zeta_n} | \hat V^{-1}_G(\btheta)-V^{-1}_G(\btheta)| =o_p(1),$ and
$\sup_{|\btheta-\btheta_0| \leq \zeta_n} | V^{-1}_G(\btheta)-V^{-1}_G(\btheta_0)|=o_p(1). $ Similar to the proof in Corollary 1, we will have  $ \hat V_G(\btheta)= V_G(\btheta)+o_p(1)$ for all $\btheta$. Assuming $\| V_G(\btheta)\|$ is bounded away from 0, so
$
   \hat V^{-1}_G(\btheta)=\lk \bI+ V^{-1}_G(\btheta)\left( \hat V_G(\btheta)-V_G(\btheta) \right) \rk^{-1}V_G^{-1}(\btheta)=O_p(1)
$
and $
    \sup_{| \btheta-\btheta_0| \leq \zeta_n} $ $| \hat V^{-1}_G(\btheta)-V^{-1}_G(\btheta)| = \sup_{| \btheta-\btheta_0| \leq \zeta_n} |\hat V^{-1}_G(\btheta)\left[\hat V_G(\btheta)-V_G(\btheta)\right]V^{-1}_G(\btheta)|=o_p(1).
$    
By Taylor expansion,
 \begin{eqnarray*}
 &&V_G(\btheta)-V_G(\btheta_0)\\
&=& E\lk\xi_i^{\otimes2}(\btheta)-\xi_i^{\otimes2}(\btheta_0)  \rk - \left\{E^{\otimes2}\lk \xi_i(\btheta)\rk -E^{\otimes2}\lk\xi_i(\btheta_0)  \rk \right\}\\
&=&E\lk 2\xi_i(\tilde\btheta)\dot \xi_{i,\theta}(\tilde\btheta)\rk (\btheta-\btheta_0)  -2E\lk \xi_i(\tilde\btheta)\rk E\lk \dot \xi_{i,\theta}(\tilde\btheta) \rk(\btheta-\btheta_0)=o_p(1),
 \end{eqnarray*}  
where $\tilde \btheta$ lies between $\btheta$ and $\btheta_0$, and $\dot \xi_{i,\theta}(\btheta)=\frac{\partial \xi_i(\btheta)}{\partial\btheta} =\dot\bg_\theta(y_i,x_i;\btheta)+(1-\delta_i) \left\{  E\lk \dot\bg_\theta(y_i,x_i;\btheta)\rk-\dot \bg_\theta(y_i,x_i;\btheta)\right\} +\delta_i C_p \frac{\partial h_n(x_i,y_i,\btheta)}{\partial \btheta}\bB(x_i)$. Thus we have
$ 
   \sup_{| \btheta-\btheta_0| \leq \zeta_n} | V^{-1}_G(\btheta)-V^{-1}_G(\btheta_0)|=\sup_{|\btheta-\btheta_0| \leq \zeta_n} | V_G(\btheta)\left[V_G^{-1}(\btheta)-V_G(\btheta_0) \right]$ $ V^{-1}_G(\btheta_0) | = o_p(1). 
$

\pagebreak

\begin{table}\caption{\label{table1} The Monte Carlo relative biases and variances of the seven estimators for the \textit{linear} and \textit{bump} models.  The number of replicates in the Monte Carlo is 1000 and the sample size is 200. $J$ is the number of imputed values. }
\tabcolsep 5.8pt

\begin{center}

(a). Model \textit{linear}: $m(x)=1+2(x-0.5)$

\vspace{0.3in}

\begin{tabular}{@{}cc || cc| cc | cc@{}}
  \hline\hline
 & & $\mu_y$ & & $\sigma_y$ & & $\rho$ & \\\hline

 & & RBias  & Var  & RBias & Var &RBias & Var  \\
 & & ($\times 100$) &($\times 100$) &($\times 100$)&($\times 100$)&($\times 100$)&($\times 100$) \\
 \hline
    & Full & 0.251 & 0.116 & -0.200 & 0.041 & -0.138 & 0.135 \\ 
   & Resp & 3.325 & 0.163 & -0.598 & 0.058 & -0.151 & 0.146 \\ 
   \hline
   J=10 & SQRI-GMM & 0.286 & 0.128 & -0.407 & 0.047 & -0.204 & 0.156 \\ 
   & MI & 0.254 & 0.119 & -0.197 & 0.044 & -0.134 & 0.148 \\ 
   & PFI & 0.251 & 0.120 & -0.463 & 0.044 & -0.332 & 0.149 \\ 
  & NPI-EL & 1.110 & 0.133 & -1.338 & 0.048 & -2.486 & 0.160 \\ 
   & HDFI & 0.365 & 0.121 & -1.142 & 0.046 & -1.538 & 0.144 \\ 
 \hline
 J=100&  SQRI-GMM & 0.256 & 0.121 & -0.404 & 0.045 & -0.196 & 0.152 \\ 
   & MI & 0.241 & 0.119 & -0.186 & 0.044 & -0.132 & 0.144 \\ 
   & PFI & 0.241 & 0.119 & -0.433 & 0.043 & -0.372 & 0.147 \\ 
   & NPI-EL & 1.114 & 0.140 & -1.305 & 0.047 & -2.286 & 0.151 \\ 
   & HDFI & 0.364 & 0.121 & -1.140 & 0.046 & -0.585 & 0.146 \\ 
    \hline\hline
\end{tabular}

\vspace{0.5in}

(b). Model \textit{bump}: $m(x)=1+2(x-0.5)+\exp\{-30(x-0.5)^2 \}$

\vspace{0.3in}

\begin{tabular}{@{}cc || cc| cc | cc@{}}
  \hline\hline
 & & $\mu_y$ & & $\sigma_y$ & & $\rho$ & \\\hline
& & RBias & Var & RBias & Var &RBias & Var  \\
 & &($\times 100$)&($\times 100$)&($\times 100$)&($\times 100$)&($\times 100$)&($\times 100$)\\
 \hline
& Full & 0.027 & 0.056 & -0.527 & 0.079 & -0.635 & 0.213 \\ 
   & Resp & 0.906 & 0.071 & -3.934 & 0.106 & -1.449 & 0.312 \\
   \hline 
 J=10  & SQRI-GMM & 0.033 & 0.064 & -0.790 & 0.083 & -0.604 & 0.224 \\ 
   & MI & 0.072 & 0.071 & -3.814 & 0.100 & -1.417 & 0.293 \\ 
   & PFI & 0.084 & 0.072 & -4.176 & 0.099 & -1.424 & 0.284 \\ 
   & NPI-EL & 0.314 & 0.061 & -3.542 & 0.093 & -4.712 & 0.270 \\ 
   & HDFI & 0.244 & 0.062 & -3.555 & 0.097 & -2.317 & 0.254 \\ 
  \hline 
 J=100  &  SQRI-GMM & 0.018 & 0.059 & -0.768 & 0.082 & -0.619 & 0.224 \\ 
   & MI & 0.077 & 0.070 & -3.833 & 0.099 & -1.355 & 0.281 \\ 
   & PFI & 0.084 & 0.070 & -4.150 & 0.099 & -1.412 & 0.280 \\ 
   & NPI-EL & 0.316 & 0.061 & -3.492 & 0.091 & -4.689 & 0.265 \\ 
   & HDFI & 0.239 & 0.061 & -3.528 & 0.096 & -2.358 & 0.254 \\  
  \hline\hline
\end{tabular}

\end{center}
\end{table}

\pagebreak

\begin{table}\caption{\label{table2} The Monte Carlo relative biases and variances of the seven estimators for the \textit{cycle} and \textit{bivariate} models.  The number of replicates in the Monte Carlo is 1000 and the sample size is 200. $J$ is the number of imputed values.}
\tabcolsep 5.8pt

\begin{center}

(c). Model \textit{cycle}: $m(x)=0.5+2x+sin(3\pi x) $

\vspace{0.3in}

\begin{tabular}{@{}cc || cc| cc | cc@{}}
  \hline\hline
 & & $\mu_y$ & & $\sigma_y$ & & $\rho$ & \\\hline
 & & RBias & Var & RBias & Var &RBias & Var  \\
 & &($\times 100$)&($\times 100$)&($\times 100$)&($\times 100$)&($\times 100$)&($\times 100$)\\ \hline
  &  Full & 0.037 & 0.182 & 0.092 & 0.065 & -0.025 & 0.047 \\ 
&    Resp & 1.942 & 0.266 & 1.973 & 0.082 & 1.177 & 0.057 \\ 
\hline
J=10  &  SQRI-GMM & 0.057 & 0.197 & -0.193 & 0.066 & -0.024 & 0.050 \\ 
   & MI & -0.200 & 0.211 & 2.500 & 0.086 & 1.373 & 0.058 \\ 
   & PFI & -0.211 & 0.210 & 2.207 & 0.086 & 1.250 & 0.058 \\ 
   & NPI-EL & 0.115 & 0.198 & -0.592 & 0.073 & -1.773 & 0.060 \\ 
   & HDFI & 0.219 & 0.187 & -0.604 & 0.070 & -1.054 & 0.056 \\  \hline
J=100  & SQRI-GMM & 0.058 & 0.185 & -0.172 & 0.066 & -0.037 & 0.050 \\ 
   & MI & -0.220 & 0.208 & 2.508 & 0.083 & 1.387 & 0.056 \\ 
   & PFI & -0.215 & 0.209 & 2.190 & 0.083 & 1.218 & 0.056 \\ 
   & NPI-EL & 0.111 & 0.194 & -0.612 & 0.073 & -1.779 & 0.059 \\ 
  & HDFI & 0.222 & 0.187 & -0.608 & 0.070 & -1.096 & 0.056 \\ 

  \hline\hline
\end{tabular}

\vspace{0.5in}

(d). Model \textit{bivariate}: $m(x)=1+2(x_1-0.5)+2\exp\{-10(x_2-0.4)^2 \}$

\vspace{0.3in}

\begin{tabular}{@{}cc || cc| cc | cc|cc@{}}
  \hline\hline
 & & $\mu_y$ & & $\sigma_y$ & & $\rho_1$ & & $\rho_2$ &\\\hline
 & & RBias & Var  & RBias & Var  &RBias  & Var & RBias  & Var  \\
 & &($\times 100$)&($\times 100$)&($\times 100$)&($\times 100$)&($\times 100$)&($\times 100$)&($\times 100$)&($\times 100$)\\
 \hline
  & Full & -0.084 & 0.305 & -0.132 & 0.100 & -0.265 & 0.152 & -0.645 & 0.453 \\ 
   & Resp & 0.676 & 0.401 & 0.186 & 0.136 & -0.724 & 0.204 & 4.101 & 0.574 \\ 
   \hline
 J=10  & SQRI-GMM & -0.094 & 0.308 & -0.369 & 0.104 & -0.295 & 0.154 & -0.704 & 0.459 \\ 
   & MI & 0.027 & 0.361 & 0.346 & 0.128 & -0.925 & 0.196 & 3.642 & 0.548 \\ 
   & PFI & 0.027 & 0.364 & 0.074 & 0.129 & -0.996 & 0.196 & 3.543 & 0.541 \\ 
   & NPI-EL & 0.306 & 0.326 & -1.090 & 0.113 & -2.177 & 0.173 & 0.562 & 0.484 \\ 
   & HDFI & 0.562 & 0.332 & -1.684 & 0.117 & -2.712 & 0.184 & 2.602 & 0.496 \\ 
    \hline
  J=100  & SQRI-GMM & -0.088 & 0.308 & -0.361 & 0.103 & -0.300 & 0.154 & -0.691 & 0.458 \\ 
   & MI & 0.024 & 0.359 & 0.346 & 0.124 & -0.898 & 0.194 & 3.615 & 0.540 \\ 
   & PFI & 0.030 & 0.358 & -0.038 & 0.125 & -0.994 & 0.194 & 3.545 & 0.538 \\ 
   & NPI-EL & 0.299 & 0.322 & -1.073 & 0.112 & -2.171 & 0.172 & 0.562 & 0.480\\ 
   & HDFI & 0.562 & 0.330 & -1.684 & 0.117 & -2.749 & 0.182 & 2.557 & 0.494 \\ 

  \hline\hline
\end{tabular}

\end{center}
\end{table}

\pagebreak

\begin{table}\caption{\label{table3} The coverage probabilities of the 95\% C.I. of the SQRI-GMM estimator under the four models. }
\begin{center}
(a). Model \textit{linear}: $m(x)=1+2(x-0.5)$

\vspace{0.1in}
\begin{tabular}{c | ccc | ccc}
  \hline\hline
  & & J=10 & & & J=100& \\
  & $\mu_y$ & $\sigma_y$ & $\rho$&$\mu_y$ & $\sigma_y$ & $\rho$ \\\hline
Normality & 0.934 &0.937 &0.817&   0.931 &0.938 &0.856 \\
Bootstrapping &    0.928& 0.953 &0.933 & 0.932 &0.953 &0.967 \\   
\hline\hline
\end{tabular}
\vspace{0.2in}

(b).  Model \textit{bump}: $m(x)=1+2(x-0.5)+\exp\{-30(x-0.5)^2 \}$
\vspace{0.1in}

\begin{tabular}{c | ccc | ccc}
  \hline\hline
  & & J=10 & & & J=100& \\
 &$\mu_y$ & $\sigma_y$ & $\rho$&$\mu_y$ & $\sigma_y$ & $\rho$ \\\hline
Normality &0.930& 0.940 &0.940&   0.947& 0.942& 0.942\\
Bootstrapping &   0.944 &0.949 &0.950 &   0.937 &0.946 &0.949\\
\hline\hline
\end{tabular}
\vspace{0.2in}

(c).  Model \textit{cycle}: $m(x)=0.5+2x+sin(3\pi x)$

\vspace{0.1in}
\begin{tabular}{c | ccc | ccc}
  \hline\hline
  & & J=10 & & & J=100& \\
  & $\mu_y$ & $\sigma_y$ & $\rho$&$\mu_y$ & $\sigma_y$ & $\rho$ \\\hline
Normality &  0.944& 0.939& 0.913&  0.943& 0.941 &0.915 \\
Bootstrapping &  0.943& 0.947& 0.941& 0.932& 0.947 &0.943 \\ 
\hline  \hline
\end{tabular}
\vspace{0.2in}

(d). Model \textit{bivariate}: $m(x)=1+2(x_1-0.5)+2\exp\{-10(x_2-0.4)^2 \}$

\vspace{0.1in}
\begin{tabular}{c | cccc | cccc}
  \hline\hline
  & & J=10 && & & J=100&& \\
  & $\mu_y$ & $\sigma_y$ & $\rho_1$&$\rho_2$&$\mu_y$ & $\sigma_y$& $\rho_1$&$\rho_2$ \\\hline
Normality&  0.953 & 0.923 &0.972& 0.939 &   0.953 &0.928 &0.977 &0.942 \\
Bootstrapping &   0.953 &0.945 &0.963& 0.950 & 0.948& 0.944& 0.958& 0.947 \\  
\hline  \hline
\end{tabular}
\end{center}
\end{table}

\begin{table}\caption{\label{table4} Relative biases and 95\% C.I. widths for the five imputation estimators in the case study. The relative biases is defined as $(\hat\btheta_n-\hat\btheta_0)/\hat\btheta_0$, where $\hat\btheta_0$ is the estimate based on full observations. The number of imputed values is $J=100$. }
\begin{center}
\begin{tabular}{c||ccc|ccc|ccc}
  \hline\hline
  &$\mu_y$& & & $\sigma_y$&& &$\rho$\\
  \hline
   & Est &RBias& Width & Est &RBias& Width & Est &RBias& Width  \\ 
   &          &  $(\times 100)$&  &&  $(\times 100)$& & & $(\times 100)$&\\
   \hline
   Full $(\hat\btheta_0)$ &  13.49 & &&0.636 & & & 0.231 &   \\
   SQRI-GMM & 13.46 &0.22 &0.15 & 0.630  &0.95  & 0.153 &0.242& 4.75&0.362\\
    MI&  13.48 &0.07 & 0.20 &  0.623&2.01 & 0.181& 0.345& 49.24&0.301 \\
    PFI& 13.48 &0.07&0.20 &   0.614&3.55 & 0.201& 0.331 & 42.83&0.371\\
    NPI-EL &13.49 &0.00&0.19 &0.594&6.59 &0.179& 0.296& 27.85&0.396\\
    HDFI&  13.49& 0.00&0.31&  0.595&8.66& 0.292&  0.306&32.11&0.503\\  
  \hline\hline
\end{tabular}
\end{center}
\end{table}


\begin{figure} \caption{\label{fig1} The comparisons of relative biases under the four models. The y-axis is for the absolute ratio between relative biases of other estimators and that of the SQRI-GMM estimator, and the x-axis is for different parameters. Curves over the horizontal line of 1 indicate the superiority of the SQRI-GMM estimator in relative biases. }
\begin{center}
\includegraphics[scale=0.34]{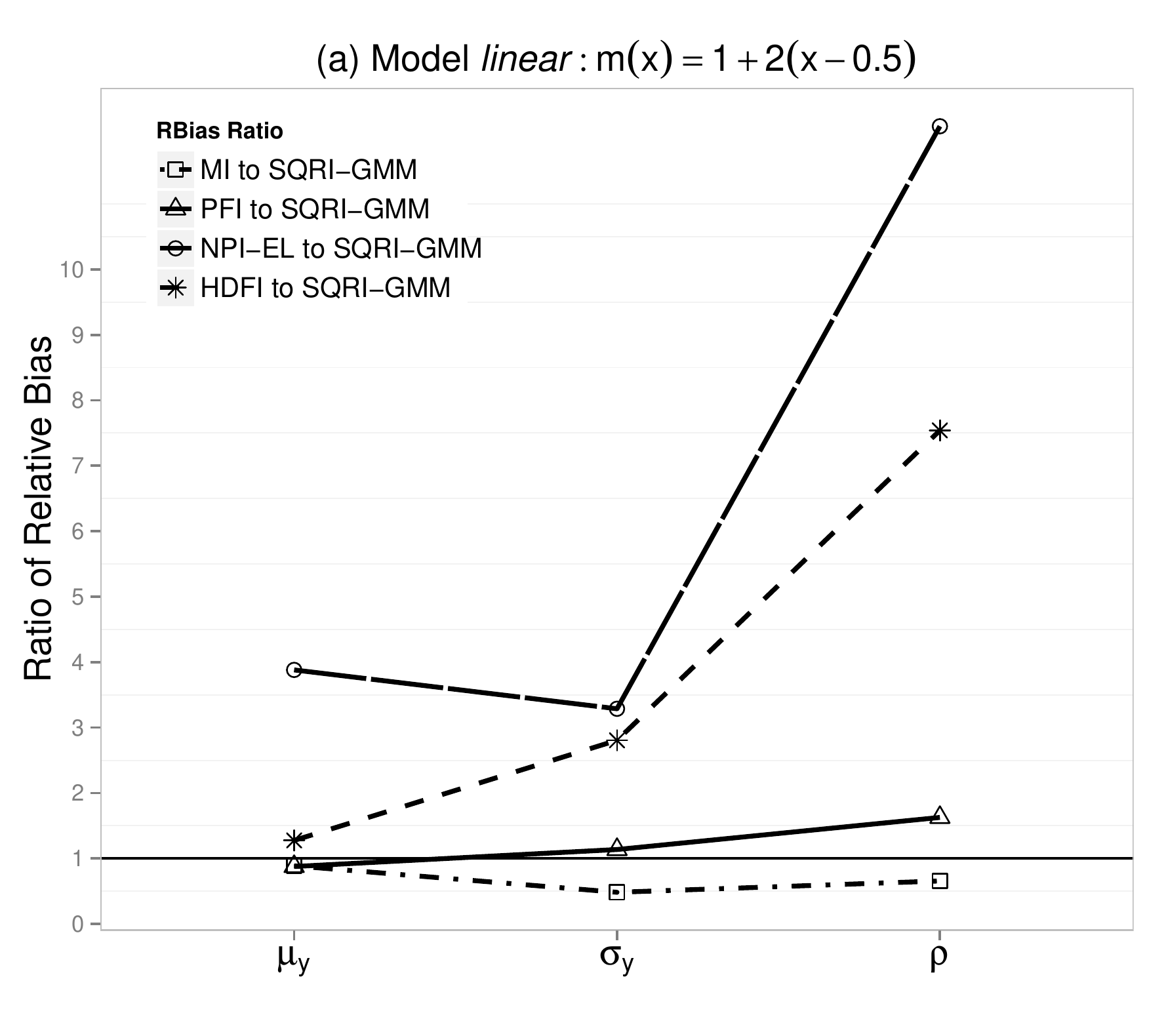}
\includegraphics[scale=0.34]{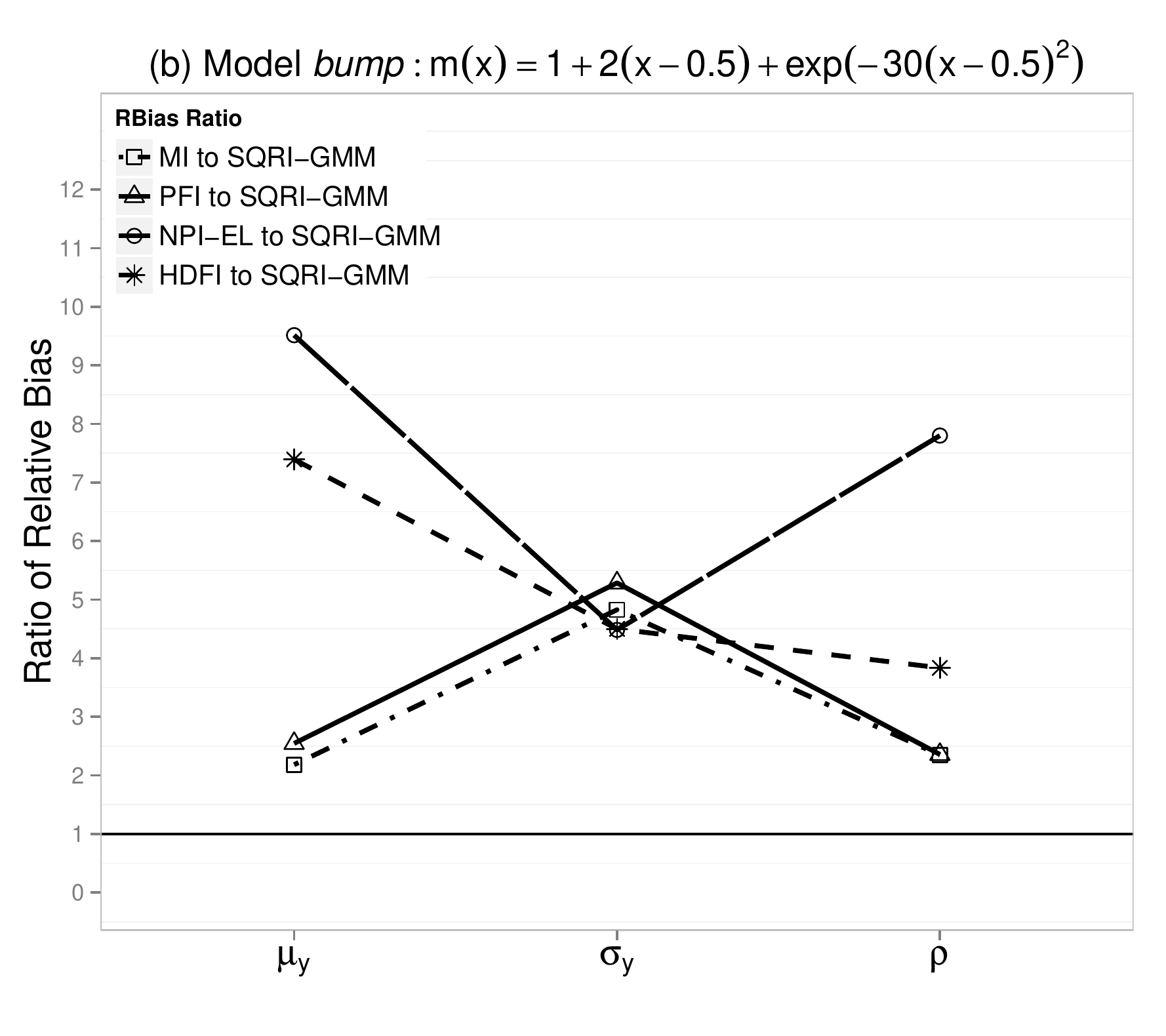}

\includegraphics[scale=0.34]{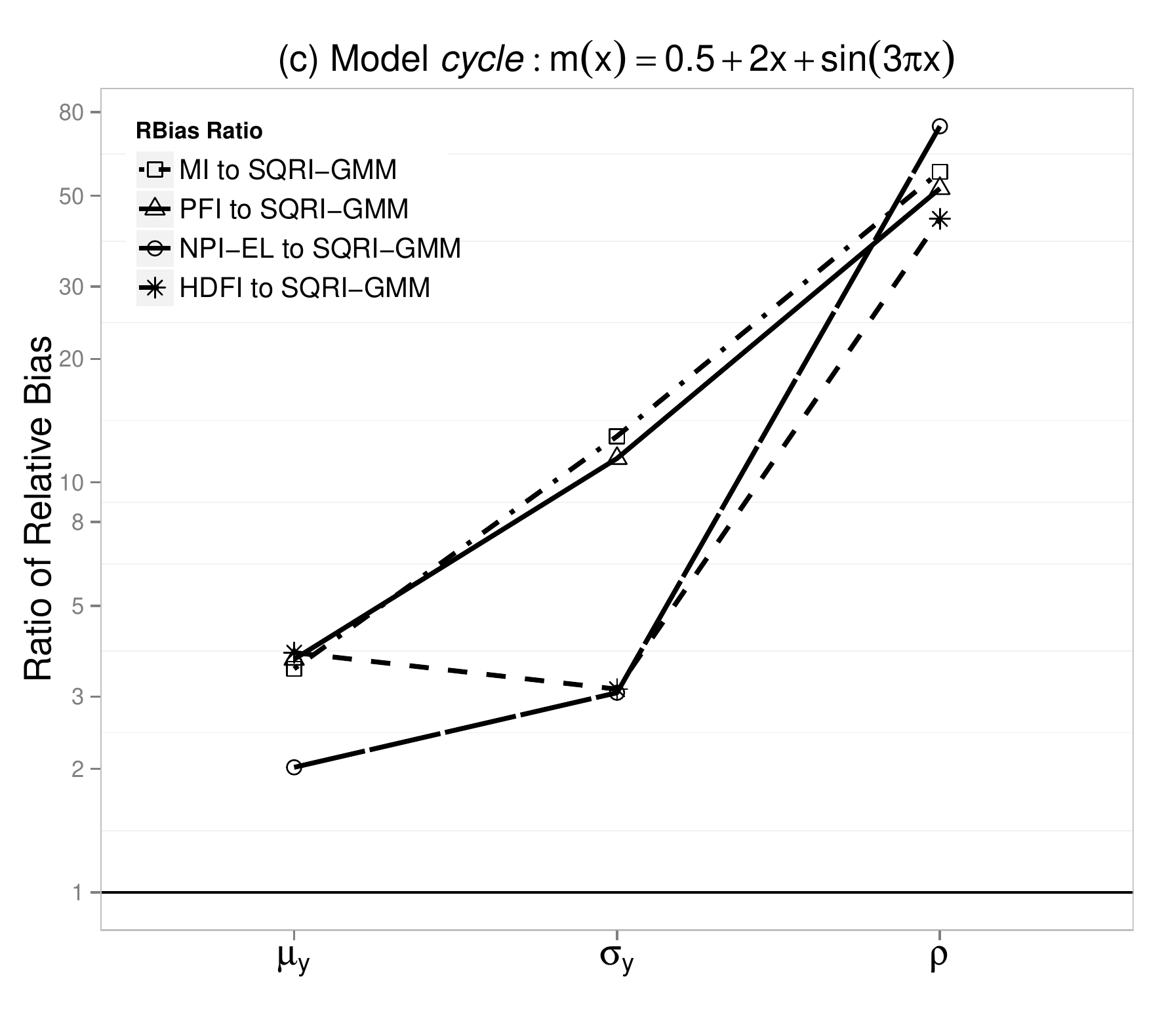}
\includegraphics[scale=0.34]{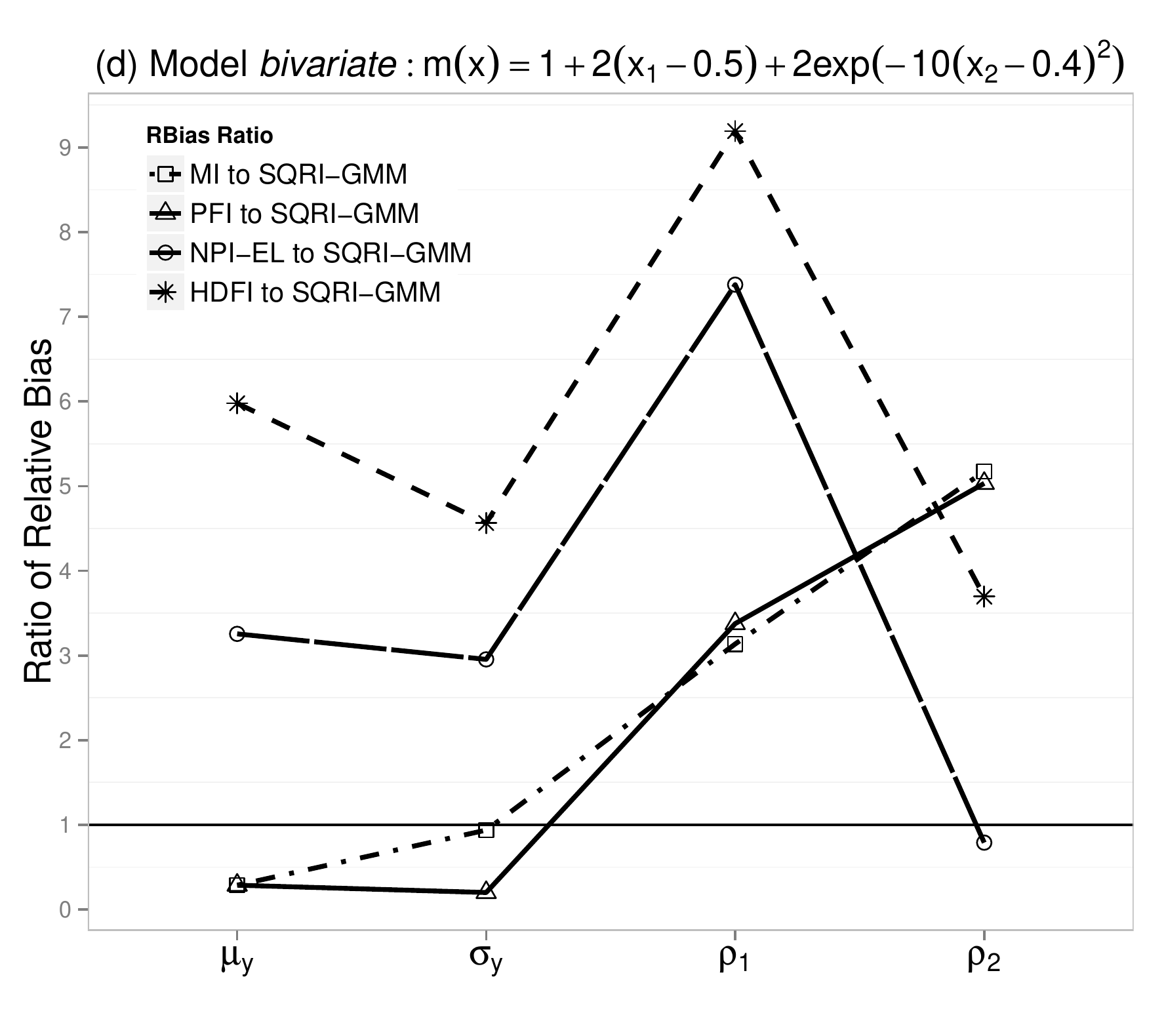}

\end{center}
\end{figure}

\begin{figure} \caption{\label{fig2} An artificial example to explain the finite sample biases observed under the two non-parametric imputation methods.}
\begin{center}
\includegraphics[scale=0.45]{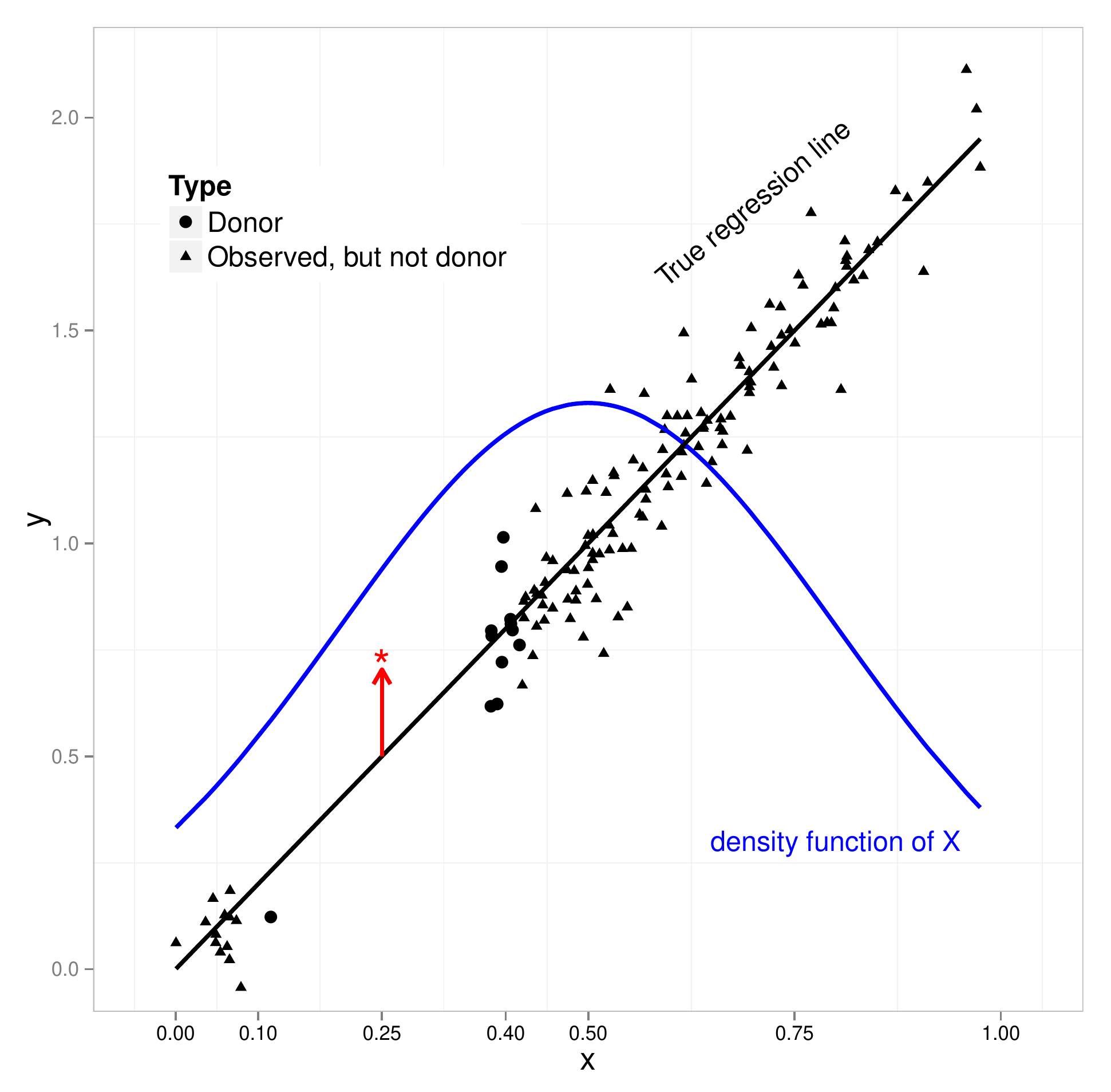}
\end{center}
\end{figure}
\begin{figure} 
\begin{center}
\caption{\label{fig3} The scatterplot of $log$(income) versus age in the case study. }

\includegraphics[scale=0.45]{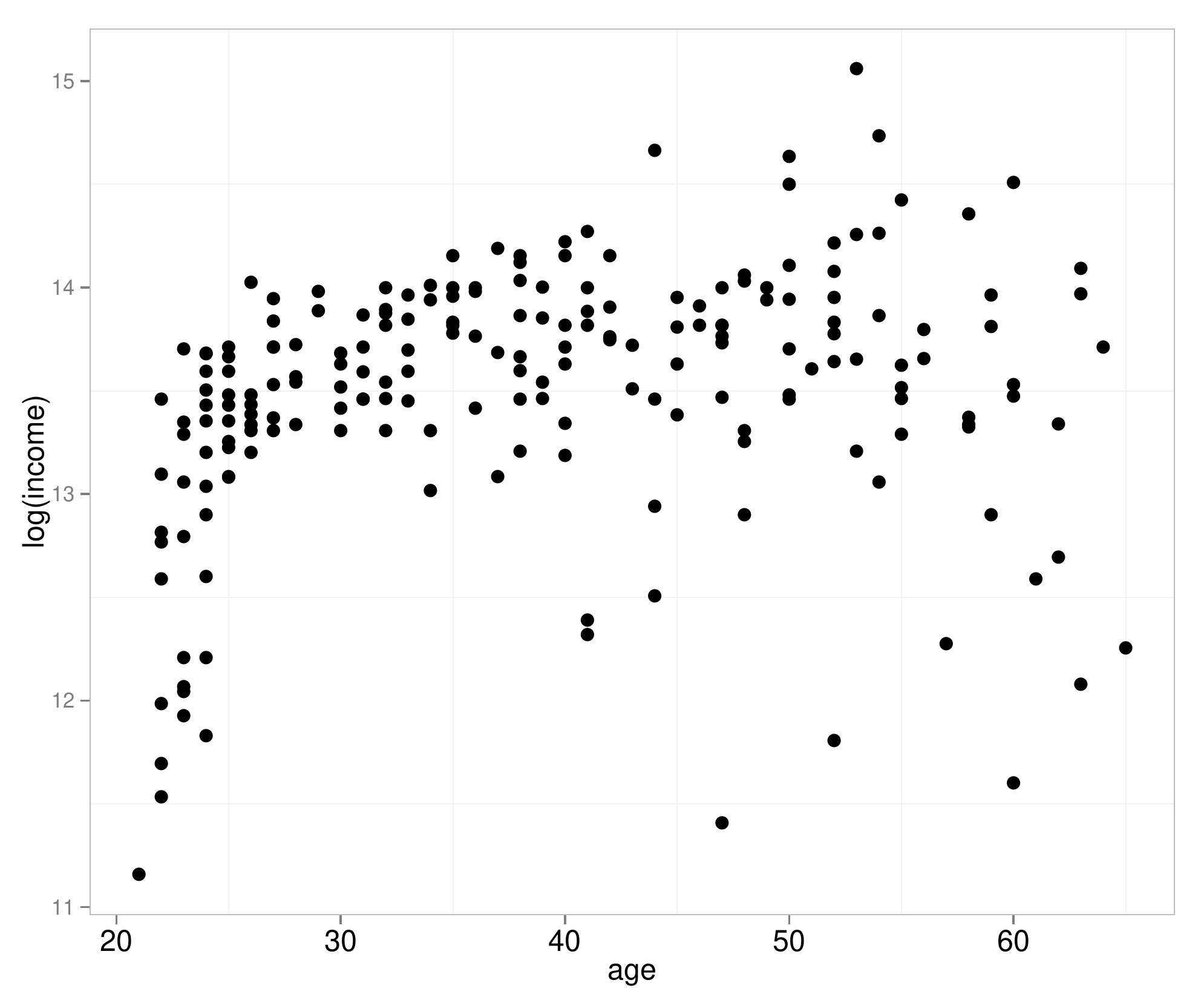}

\end{center}
\end{figure}
\end{document}